\documentclass[a4paper,10pt]{article}

\usepackage{graphicx}
\usepackage{amsfonts, amssymb, amsmath, fullpage}

\usepackage{tikz}

\title{Stability of Periodic Travelling Flexural-Gravity Waves in Two Dimensions}

\author{Olga Trichtchenko, Paul Milewski, Emilian P\u{a}r\u{a}u and Jean-Marc Vanden-Broeck}
\date{}

\begin{document}

\maketitle

\abstract{In this work, we solve  the Euler's equations for periodic waves  travelling under a sheet of ice using a reformulation introduced in \cite{AFM06}.  These waves are referred to as flexural-gravity waves. We compare and contrast two models for the effect of the ice:  a linear model and a nonlinear model. The benefit of this reformulation is that it facilitates the asymptotic analysis. We use it to  derive the nonlinear Schr\"{o}dinger equation that describes the modulational instability of periodic travelling waves. We compare this asymptotic result with numerical computation of stability using the Fourier-Floquet-Hill method and show how well these agree. We show that different models have different stability regimes for large  values of the flexural rigidity parameter. Numerical computations are used to go beyond the modulational instability and show high frequency instabilities that are the same for both models for ice in the regime examined.}

\section{Introduction}\label{sec:introduction}
In this work, we examine water waves under a sheet of ice, referred to as flexural-gravity waves or hydroelastic waves. We model the water as an incompressible, inviscid and irrotational fluid, restricting our focus on two-dimensional waves with periodic boundary conditions. The contribution from the ice to the movement of the wave can be modelled in several ways with some models that conserve energy and some that do not. One of the earliest instances of modelling ice was perhaps shown in the paper by Greenhill \cite{G86}  (for a more complete review, see Squire et al. \cite{SDWRL95, S07}). The two Hamiltonian (conservative) models \cite{MW13} we consider are the linear (biharmonic) Euler-Bernoulli model and the nonlinear model derived from the Cosserat theory of shells \cite{PT11}, which can also be considered as a Willmore functional \cite{W82} using the formulation by Toland \cite{T08}. In this work, we compare and contrast solutions to these models and the stability of these solutions. We are interested in the shape of the interface, which makes this a free boundary problem. Furthermore, we restrict the problem to waves moving at a constant speed. 

There are traditionally two categories of waves studied, periodic as in this work, and solitary waves on an infinite domain. Solitary waves can either be forced, for example by a moving load on top of the ice, or free solitary waves. These can further be split into two regimes, deep water and finite depth water.  Once a model is proposed, a natural question to ask is  whether or not the model admits solutions. Several works discuss the  existence of solutions to equations describing hydroelastic waves; for example the work of Toland \cite{T08} discusses the  existence of solutions as an optimization of the Lagrangian formulation for travelling waves. Using a variational approach, Groves et al \cite{GHW16} show the existence for hydroelastic solitary waves and Akers et al. \cite{AAS17} use bifurcation theory for  the existence of periodic waves in two dimensions.

Without solving the full set of the proposed equations, some  insight  can be gained from asymptotic models for small amplitude solutions. Using more physical considerations and supplementing with observational results, Liu and Mollo-Christensen \cite{LMC88} derive a weakly nonlinear form of the governing equations for waves in an ice pack, including their stability analysis. Using Hamiltonian formalism, Marchenko and Shrira \cite{MS91} consider weakly nonlinear theory and determine the stresses in the ice.
Forced waves in water of finite and infinite depth were considered by P{\u{a}}r{\u{a}}u and Dias \cite{PD02}. Using  the normal form theory and considering travelling waves moving close to a critical speed of the wave for one model for the  ice, the analysis was reduced to studying the nonlinear Schr\"{o}dinger equation to show that below a critical depth, solitary wave solutions exist. Considering three different models for the  ice including the one in \cite{PD02} and a different asymptotic reduction, Milewski and Wang \cite{MW13} studied flexural-gravity solitary waves in two and three dimensions and concluded for there were no small amplitude solutions for certain values of parameters. 

Using the boundary integral method, Vanden-Broeck and P{\u{a}}r{\u{a}}u \cite{VBP11} were  able to compute both solitary and periodic travelling waves for a simpler nonlinear model originally proposed by Forbes \cite{F86}. The fully nonlinear model was considered in Gao and Vanden-Broeck \cite{GVB14} also for both periodic and solitary waves. A more general discussion considering periodic interfacial waves with and without mass is seen in Akers et al. \cite{AAS17, AAS172} where a different parametrisation of the problem was considered. Work on computing solutions for the three-dimensional problem for flexural-gravity waves also exists, but will not be discussed here. 

The presence of ice introduces more nonlinearity and higher order derivatives than  in previous work for gravity-capillary waves \cite{DT14}, but we can  follow  a similar methodology to reformulate the problem. In the presence of the flexural term, resonance similar to that first observed by Wilton \cite{W15} for capillary-gravity waves occurs for specific parameters. They are of a similar nature for both gravity-capillary waves and flexural-gravity waves. These can be studied  numerically by the methods introduced in \cite{TDW16}. 

In this paper, we use the reformulation introduced by Ablowitz, Fokas and Musslimani \cite{AFM06} and extend it for different conditions at the surface. This is useful not only for performing an {asymptotic analysis} in the regime where the nonlinear Schr\"{o}dinger (NLS) equation applies, but also for computing more general results numerically. The NLS equation  allows us to compare the  modulational instability (derived asymptotically)  to the stability results (computed numerically) for flexural-gravity waves to see how well these match for different models. The second type of instabilities referred to as high-frequency instabilities \cite{DT16} are also examined numerically using the Fourier-Floquet-Hill method for the time dependent problem \cite{DK06} using different models to describe the ice at the surface.

The outline of the paper is as follows. In Section \ref{sec:model}, we present the different models used to describe how water waves behave  under ice and reformulate the equations into a form which facilitates asymptotic approximations in Section \ref{sec:asymptotics}. Following these analytical results, we show how these reformulated equations can be solved numerically and set up the numerical spectral stability eigenvalue problem in Section \ref{sec:numerics}. In Section \ref{sec:NLSresults} we show that in the proper limit, the numerical results agree with those from the asymptotic analysis. Numerical computations are then used more generally, giving a richer understanding of the types of solutions and of their stability in Section \ref{sec:genResults}.  We conclude in Section \ref{sec:conclusion}.

\section{Model and Reformulation} \label{sec:model}
To model irrotational, inviscid and incompressible flows  under a variety of surface conditions, we use the  Euler's equations written in potential form as
\begin{align}
\begin{cases}
\displaystyle
\phi_{xx} + \phi_{zz}  = 0,  \ \ & (x,z) \in \mathcal{D}, \\
\displaystyle
\phi_z  = 0,  \ \ &z=-h \\
\displaystyle
 \eta_t + \eta_x \phi_x  = \phi_z,  \ \ &z = \eta(x,t) \\
 \displaystyle
 \phi_t + \frac{1}{2} \left( \phi_x^2 + \phi_z^2 \right) + g \eta = -\frac{D}{\rho}P_{flex},  \ \ &z=\eta(x,t)
\end{cases}
\label{eq:euler}
\end{align}
where $h$ is the height of the fluid, $g$ is the acceleration due to gravity, $\eta(x,t)$ is the elevation of the fluid surface,  $\rho$ is the density and $\phi(x,z,t)$ is  the  velocity potential. 
Here $D$ is the flexural rigidity defined by
\begin{align*}
D = \frac{E d^3}{12(1-\nu^2)}
\end{align*}
where  $E$ is  the Young's modulus, $\nu$ is the Poisson ratio and $d$ is  the thickness of the ice.

\begin{figure}[h!]
\begin{center}
\begin{tikzpicture}
\fill[fill=gray!20, scale=0.67] (0,0.2) -- plot [domain=0:12] ({\x},{cos(deg(\x))/4+2.9}) -- (12,0.2) -- cycle;
\draw[fill=black!80, color = black](0,0) rectangle (8,0.2);
\draw[scale=0.67,domain=0:12,smooth,variable=\x]  plot ({\x},{cos(deg(\x))/4+2.9});
\draw[scale=0.67,domain=0:12,smooth,variable=\x]  plot ({\x},{cos(deg(\x))/4+3.1});
\draw[->] (0,0.2) -- (8.2,0.2);
\draw[->] (0,0.2) -- (0,3.2);
\node at (8.5,0.2) {$x$};
\node at (0,3.5) {$z$};
\node at (-0.7,2.5) {$\eta(x,z,t)$};
\node at (-0.7,0.2) {$z=-h$};
\draw[color = black, fill=white!60](3.5,0.85) rectangle (4.5,1.4);
\node at (4,1.1) {water};
\draw[color = black, fill=white!60](3.6,1.9) rectangle (4.4,2.4);
\node at (4,2.15) {ice};
\draw[color = black, fill=white!60](3.6,2.75) rectangle (4.4,3.2);
\node at (4,3) {air};
\end{tikzpicture}
\caption{A schematic of the physical scenario. \label{fig:schematic}}
\end{center}
\end{figure}
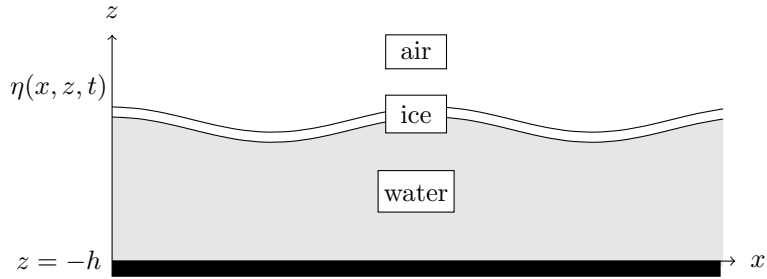

We focus on solutions on a periodic domain with the schematic shown in Figure \ref{fig:schematic} where the unknown domain $\mathcal{D}$ is shown in grey. In the Bernoulli equation (the last equation in (\ref{eq:euler})),  $P_{flex}$ can represent a variety of conditions at the surface, such as surface tension or in our case, a thin sheet of ice on the surface of the water. Several models exist that are based on considering the ice as an elastic sheet and we focus on the following: 
\begin{enumerate}
\item The linear (biharmonic) model assuming  that the ice behaves like an Euler-Bernoulli thin elastic plate in the regime where the curvature is small
\begin{align}
P_{flex} = \eta_{4x}.
\end{align}
\item The nonlinear (Toland or Cosserat) model as shown in \cite{PT11}. It is given by
\begin{align}
P_{flex} = \partial^2_x\left[ \frac{\eta_{xx}}{(1+\eta_x^2)^{5/2}}\right] + \frac{5}{2}\partial_x\left[ \frac{\eta^2_{xx}\eta_x}{(1+\eta_x^2)^{7/2}}\right].
\end{align}
\end{enumerate}
In these models, the ice  is assumed to be a thin elastic plate with constant thickness which bends with the water waves. Furthermore, the friction between the ice and the water is neglected.

We are interested in  studying how the interface $\eta(x,t)$ changes depending on the model used. We  rewrite the equations solely in terms of the surface variables \cite{Z68} by introducing the velocity potential defined at the surface {$q(x,t) = \phi(x,\eta(x,t),t)$}. Applying the  chain rules to compute all the derivatives as shown in \cite{AFM06} and making use of the third equation in \eqref{eq:euler} which is valid at the surface,  allows us to rewrite the Bernoulli equation (the last equation in \eqref{eq:euler}) as
\begin{align}
 q_t + \frac{1}{2} q_x^2 + g\eta -
\frac{1}{2}\frac{(\eta_t+\eta_xq_x)^2}{1+\eta_x^2} = - \hat{D}
P_{\text{flex}}
\label{eq:local}.
\end{align}
We refer to \eqref{eq:local} as the local equation. We introduced the notation $\hat{D} = \frac{D}{\rho}$ and will subsequently drop the hat notation for ease. There are several reformulations that exist for the given set of Euler's equations (for a review see \cite{WV14}). In this work, the focus is on the implicit formulation that we will refer to as Ablowitz-Fokas-Musslimani (AFM) formulation. In  \cite{AFM06}, the authors introduce an identity for two functions satisfying Laplace's equation and write the expression in divergence form. Choosing a particular solution for one of the functions and applying  the divergence theorem,  defines a global relation often seen in the unified transform method of Fokas \cite{F97} for solving partial differential equations. Making use of the boundary conditions as well as the  periodicity \cite{DO11}, we obtain a  nonlocal equation given by
\begin{align}
\int_0^Le^{ikx}\left( i\eta_t \cosh(k(\eta+h)) + q_x\sinh(k(\eta+h)) \right)dx = 0, ~~ \ k \in \Lambda,
\label{eq:nonlocal}
\end{align}
where $\Lambda$ is defined as the period lattice given by $\Lambda =  \{2\pi n/L \ | \ n \in \mathbb{Z}, n \ne 0\}$ and $L$ the period of the solution. We now restrict $L$ to be $2\pi$.

\section{Asymptotic Analysis} \label{sec:asymptotics}
While the system is in a more compact form, it is still hard to solve the set of time-dependent equations given by \eqref{eq:local} and \eqref{eq:nonlocal}. We examine the local equation \eqref{eq:local} and the nonlocal equation \eqref{eq:nonlocal} asymptotically. Following the work of \cite{AFM06}, we focus on small-amplitude, slowly varying envelope equation for quasi-monochromatic waves. For the purposes of this section only, we restrict this analysis to infinitely deep water $h = \infty$, where (as it will be shown), the governing equation of motion is the nonlinear Schr\"{o}dinger equation (NLS). This allows us to obtain the parameter regime for modulational instability (the focussing case of NLS as will be defined later) or midulational stability (the defocussing case). We extend the procedure outlined in \cite{AFM06} where gravity-capillary waves were considered.  In this work, a detailed derivation is presented since the authors feel that this clarification is needed.

The local equation is the same in finite and infinite depth and given by \eqref{eq:local}, whereas the nonlocal equation in infinite depth is given by
\begin{align*}
\int_{0}^{2\pi}dxe^{-ikx}e^{|k|\eta}\left[ i\eta_t - \text{sgn}(k)q_x \right] = 0.
\end{align*}
First, we assume small-amplitude solutions setting $\eta \rightarrow \epsilon \eta$ and $q \rightarrow \epsilon q$. Keeping terms of up to second order in the small parameter $\epsilon$, we obtain 
\begin{align*}
\int_{0}^{2\pi}dxe^{-ikx}\left[ 1 + \epsilon|k|\eta + \epsilon^2\frac{1}{2}|k|^2\eta^2\right]\left( i\eta_t - \text{sgn}(k)q_x\right) & = 0 \\
q_t + g\eta +DP_{\text{flex}}(\epsilon \eta) + \frac{1}{2}\epsilon q_x^2 - \frac{\epsilon}{2}(\eta_t^2 +2\epsilon\eta_t\eta_xq_x) & = 0
\end{align*}
where, for example, the linear model is
\begin{align*}
P_{\text{flex}}(\epsilon \eta) = \eta_{4x}.
\end{align*}

To simplify, we differentiate the local equation with respect to $x$ and let $Q = q_x$
\begin{align*}
\int_{0}^{2\pi}dx e^{-ikx}\left( i \eta_t - \text{sgn}(k) Q + i \epsilon |k| \eta\eta_t - \epsilon k \eta Q + \frac{1}{2} i k^2 \epsilon^2 \eta^2 \eta_t  - \frac{1}{2} \epsilon^2 |k| k \eta^2 Q \right) & = 0 \\
Q_t + g\eta_x + D P_{x,\text{flex}}(\epsilon \eta) + \epsilon QQ_x - \epsilon \eta_t\eta_{tx} - \epsilon^2 \eta_{tx}\eta_xQ - \epsilon^2\eta_t\eta_{xx}Q - \epsilon^2\eta_t\eta_x Q_x & = 0.
\end{align*}
Focussing on waves with slow varying envelopes  and rapidly oscillating carrier waves, we can now introduce slow and fast variables $X = \epsilon x$ and $T = \epsilon t$ such that $\partial_x \rightarrow \partial_x + \epsilon \partial_X$ and $\partial_t  \rightarrow \partial_t + \epsilon \partial_T$ to obtain up to $O(\epsilon^2)$
\begin{align*}
\int_{0}^{2\pi}dxe^{-ikx}\left( i\eta_t  - \text{sgn}(k)Q + \epsilon i\eta_T + \epsilon i |k|\eta \eta_t  - \epsilon k\eta Q + \epsilon^2 i |k|\eta \eta_T + \epsilon^2\frac{i}{2}k^2\eta^2 \eta_t - \epsilon^2\frac{1}{2}\text{sgn}(k)k^2\eta^2 Q\right) & = 0 \\
Q_t + g\eta_x  + DP_{x,\text{flex}}(\epsilon \eta) + \epsilon g \eta_X + \epsilon Q_T +  \epsilon QQ_x - \epsilon \eta_t \eta_{tx}  - \epsilon^2 \eta_T \eta_{tx} -\epsilon^2 \eta_t \eta_{Xt}  - \epsilon^2 \eta_t \eta_{xT} + \epsilon^2 QQ_X   \\ - \epsilon^2 \eta_{tx}\eta_xQ - \epsilon^2\eta_t\eta_{xx}Q - \epsilon^2\eta_t\eta_xQ_x & = 0
\end{align*}
with, for example, the derivative of the biharmonic term given by
\begin{align*}
P_{x,\text{flex}}(\epsilon \eta) =  \eta_{5x} + 5 \epsilon \eta_{4xX} + 10 \epsilon^2 \eta_{3x2X}.
\end{align*}
We  note that in \cite{AFM06}  the approximation to the local equation is missing a bracket.
We now look for solutions in terms of quasi-monochromatic waves of the form
\begin{align}
\eta & = \eta_1e^{i\theta}  + \epsilon\left(\eta_0 + \eta_2e^{2i\theta}\right) + O(\epsilon^2) + c.c.\label{eq:waveDecomp}\\
Q & = Q_1e^{i\theta} + \epsilon\left(Q_0 + Q_2e^{2i\theta} \right) + O(\epsilon^2) + c.c.
\end{align}
where  $\eta_j = \eta_j(X,T)$ and $Q_j = Q_j(X,T)$ and $\theta = kx-\omega t$ and c.c. denotes the complex conjugate. Using the above it is important to note that up to first order
\begin{align*}
\partial_T\eta & = \eta_{1,T}e^{i\theta} +  \bar{\eta}_{1,T}e^{-i\theta} + \epsilon\left(\eta_{0,T} + \eta_{2,T}e^{2i\theta} +  \bar{\eta}_{2,T}e^{-2i\theta}\right) \\
\partial_X\eta & = \eta_{1,X}e^{i\theta} +  \bar{\eta}_{1,X}e^{-i\theta} + \epsilon\left(\eta_{0,X} + \eta_{2,X}e^{2i\theta} +  \bar{\eta}_{2,X}e^{-2i\theta}\right) \\
\partial_t\eta & = -i\omega \eta_{1}e^{i\theta} + i\omega \bar{\eta}_{1}e^{-i\theta} + \epsilon\left( -2i\omega \eta_{2}e^{2i\theta} + 2i \omega \bar{\eta}_{2}e^{-2i\theta}\right) \\
\partial_x\eta & = ik \eta_{1}e^{i\theta} - ik \bar{\eta}_{1}e^{-i\theta} + \epsilon\left( 2ik \eta_{2}e^{2i\theta} - 2i k \bar{\eta}_{2}e^{-2i\theta}\right),
\end{align*}
and similarly for $Q(x,t,X,T)$ where the barred quantities are the complex conjugates.

We derive an  equation for the leading terms of the wave profile $\eta(x,t,X,T)$ by grouping terms of different orders in $\epsilon$ and wavenumbers $k$.  We
now outline the  procedure,  using the simplest (linear/biharmonic) model for flexural-gravity waves as an example.
The constant terms show that $\eta_0$ and $Q_0$ are zero at lowest orders in $\epsilon$.  For the
 leading order terms ($O(\epsilon^0)$) of the coefficient for $e^{i\theta}$, we obtain
\begin{align}
-\text{sgn}(k) Q_1 + \omega \eta_1 & = 0 \\
 - i Q_{1} \omega + ikg\eta_1 + D P_{x,\text{flex}}(\eta_1) & = 0 
\end{align}
with the biharmonic term for flexural-gravity waves given by
\[P_{x,\text{flex}}(\eta_1) = i k^5 \eta_1. \]
This allows us to compute the first correction term as
\begin{align}
\omega^2 & =  \text{sgn}(k) k (k^4 D +g) \\
Q_1 & = \frac{ k(k^4 D  +g)}{\omega}\eta_1 + O(\epsilon)
\label{eq:q11}.
\end{align}
The first two terms of the coefficient of $e^{i\theta}$ in the local equation give
\begin{align}
\epsilon \left(g \eta_{1, X} + Q_{1,T} \right) - i Q_1\omega + i k g \eta_1 + D P_{x,\text{flex}} = 0, 
\label{eq:q12}
\end{align}
with \[P_{x,\text{flex}} = 5 k^4 \epsilon \eta_{1,X} + i k^5\eta_1.\]
Using \eqref{eq:q11} to substitute into the terms multiplied by $\epsilon$ (\textit{i.e} all the derivatives) we can solve for the first order correction
\begin{align}
Q_{1} = \left[Dk^5+ gk - \epsilon\frac{i}{\omega} \left( 5 k^4D\omega \frac{\partial}{\partial x} + g\omega \frac{\partial}{\partial x} + k^5 \text{sgn}(k)D \frac{\partial}{\partial t} + gk\text{sgn}(k)\frac{\partial}{\partial t} \right)\right] \frac{\eta_1}{\omega} + O(\epsilon^2).
\label{eq:q12}
\end{align}
To obtain the second order ($O(\epsilon^2)$) correction for $Q_1(X,T)$, we need $Q_2(X,T)$ and $\eta_2(X,T)$ which occurs at the highest order in $\epsilon$ in the coefficient of $e^{2i\theta}$. From the nonlocal and local equations, we obtain
\begin{align}
-\text{sgn}(k) Q_2 - k \eta_1 Q_1 + |k|\omega \eta_1^2 + 2\omega \eta_2 & = 0\\
ik Q_1^2 + DP_{x,\text{flex}}(\eta_2) + i \omega^2 k \eta_1^2 + 2ikg\eta_2 - 2i\omega Q_2  & = 0,
\end{align}
where for the biharmonic model, 
\[ P_{x,\text{flex}}(\eta_2) = 32i k^5\eta_2 .\]
From the above two equations, we obtain 
\begin{align}
\eta_2 & = \frac{g+k^4D}{g -14k^4D} |k| \eta_1^2+ O(\epsilon) \\
Q_2 & = \frac{g+k^4D}{g-14k^4D} 2\omega k \eta_1^2  + O(\epsilon).
\end{align}
Repeating the above, but for the first three terms of the coefficient of $e^{i\theta}$, we obtain from the local equations
\begin{align}
Q_1 & = \frac{k(Dk^4+ g)}{\omega}\eta_1 - \epsilon\frac{i}{\omega} \left( (5 k^4D + g)\eta_{1,X} -  \frac{ik}{\omega^2} (k^4 D + g) \eta_{1,T} \right) \nonumber \\ & - \epsilon^2 \frac{10 k^3 D}{\omega} \eta_{1,XX} - \epsilon^2 \frac{1}{\omega} \eta_{1,TT} + \epsilon^2 \frac{(k^4D + g)k^3}{\omega} + O(\epsilon^3).
\label{eq:q13}
\end{align}
 Finally, substituting \eqref{eq:q13} into the nonlocal equation and grouping the coefficients of $\eta_1(X,T)$ we obtain the nonlinear Schr\"{o}dinger (NLS) equation for the envelope of the wave profile
\begin{align}
i(\partial_T+\omega'\partial_X)\eta_1 + \epsilon \frac{\omega''}{2}\partial_X^2\eta_1 + \epsilon M|\eta_1|^2\eta_1 = 0.
\label{eq:NLSwE}
\end{align}

\noindent The linear dispersion relation $\omega$ appearing in NLS is independent of the model used and so are its derivatives and these are given by
\begin{align}
\omega^2 & = |k|(g+k^4D) \\
\omega'& = \frac{\text{sgn}(k)(g+5k^4D)}{2\omega} \\
\omega{''} & = -\frac{\omega(g^2-30gk^4D-15(k^4D)^2)}{4k^2(g+k^4D)^2}.
\label{eq:omegapp}
\end{align}
The second derivative of dispersion \eqref{eq:omegapp} is the same as in \cite{LMC88} if inertia and compression are neglected as done in this work. However, the term multiplying the nonlinearity depends on the model.  The different quantities are 
\begin{align}
M = -\frac{\omega k^2(4g^2-27gk^4D+44(k^4D)^2)}{2(g+k^4D)(g-14k^4D)}
\label{eq:MTol}
\end{align}
for the nonlinear (Toland or Cosserat) model  and
\begin{align}
M = -\frac{\omega k^2(2g^2-11gk^4D-13(k^4D)^2)}{(g+k^4D)(g-14k^4D)}
\label{eq:MLin}
\end{align}
for the linear (biharmonic) model.

To obtain the nonlinear Schr\"odinger equation in a more standard form, we introduce the group velocity $v_g = \omega'$ and a slow time and space variables $\tau = \epsilon T$ and $\xi = X-v_gT$. Once we divide through by the small parameter $\epsilon$, the nonlinear Schr\"odinger equation becomes
 \begin{align}
i\eta_{1,\tau} + \frac{\omega''}{2} \eta_{1,\xi \xi} + M|\eta_1|^2\eta_1 = 0.
\label{eq:genNLS}
\end{align}

\begin{figure}
\begin{center}
\begin{tikzpicture}
    \node(a){\includegraphics[width=0.65\textwidth]{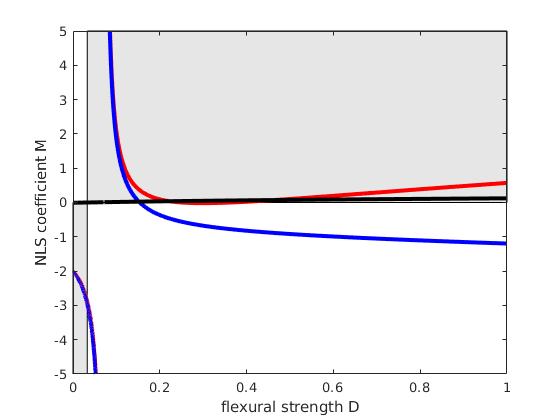}};
    \node at (2.3, 1.9)
    [anchor=center]
    {\includegraphics[width=0.2\textwidth]{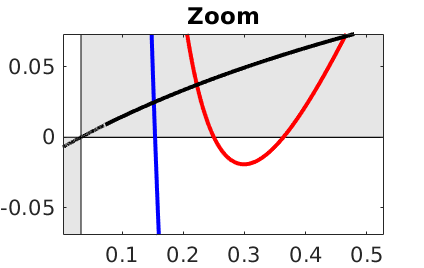}};
    \node at (1.2, -0.82) {LIN};
    \node at (1.2, 0.455) {NL};
\end{tikzpicture}
\end{center}
\caption{Coefficient $M$ in front of the nonlinear term in \eqref{eq:NLSwE}  as a function of the flexural rigidity $D$ for two different models representing flexural-gravity waves. The red is the nonlinear model (labelled NL) and the blue is the linear model (labelled LIN). The thick black line represents \ref{eq:omegapp}. The grey is the unstable region (focussing NLS regime) and the white is the stable regime.  We can see for small flexural rigidity the models go from unstable to stable in a similar way, but for large $D$, the two models differ. Inset shows interesting region near zero.
\label{fig:NLSCoeff}}
\end{figure}

The NLS equation is of focusing type when $\omega''M>0$ and also modulationally unstable \cite{MW13}. It is interesting to note that the denominator of $M$ becomes zero at $D = 1/14 \approx 0.07$ for $g=1$. In the case of gravity-capillary waves, the condition for the vanishing denominator (resonance) is referred to as a Wilton ripple  \cite{W15,TDW16} and a similar terminology will be used here. The summary of the different regions of stability and instability as determined by the coefficients of NLS for different models is shown in Figure \ref{fig:NLSCoeff}. The vertical asymptote represents the Wilton ripple.  We are interested in how varying $D$ changes the stability properties of the solutions with instability regions represented in grey. The second derivative of the dispersion, $\omega''$, changes sign once at $D \approx 0.03$ and is the change from the grey (unstable) to white (stable region), while the nonlinear coefficient $M$ stays negative. The two different lines represent the nonlinear coefficient $M$ seen in Equation \eqref{eq:genNLS} for different models with the red (labelled NL) representing the nonlinear (Toland) model and linear model (bihamornic) for ice shown in blue (labelled LIN). The greatest arises for large $D$ which represents a more rigid regime of the elastic sheet. In this case we see that the nonlinear (Toland) model is first unstable, briefly becomes stable and then goes back to being unstable for large enough $D$, whereas the linear model remains stable once it transitions.

We continue to examine how instabilities and grow in time. This can be done by noting that a spatially independent solution of \eqref{eq:genNLS} of amplitude $a$ is given by
\begin{align*}
\eta_1^{(0)}(\xi,\tau) = ae^{i M a^2 \tau},
\label{eq:NLSA0}
\end{align*}
which implies that $\eta(x,t)$ will travel at a constant speed and will be monochromatic (at first order) with wavenumber $k$ as shown in \eqref{eq:waveDecomp}. We perturb this particular solution by an arbitrary complex function of magnitude $\delta$ such that 
\begin{align}
\eta_1(\xi,\tau) = \left[a + \delta\left(f(\xi,\tau) + i g(\xi,\tau) \right)\right] e^{i M a^2 \tau},
\end{align} 
where $f$ and $g$ are arbitrary functions. Assuming the perturbation is small, then up to first order in $\delta$ we obtain the  real and imaginary parts of the perturbation as two coupled equations
\begin{align*}
g_{\tau} - \frac{\omega''}{2} f_{\xi\xi} - 2M a^2 f & = 0 \\
f_{\tau} + \frac{\omega''}{2}  g_{\xi\xi} & = 0.
\end{align*}
Since we are interested in the case where the solution to the above equation becomes unstable, we can look for the following form of the perturbation
\begin{align*}
f(\xi,\tau) & = {u}e^{\Omega \tau} e^{i\mu\xi} \\
g(\xi,\tau) & = {v}e^{\Omega \tau} e^{i\mu\xi}
\end{align*}
where $\mu$ is real, and examine when the solution will grow exponentially in time, i.e. where $\Omega$ is real and positive. The second equation gives that ${v} = \Omega/(\alpha \mu^2) {u}$ and using the first equation, we obtain the time dependence as
\begin{align}
\Omega^2 = {\omega''} M a^2 \mu ^2 - \left( \frac{\omega''}{2}\right)^2 \mu^4.
\label{eq:NLSGrowthRate}
\end{align}
The maximum of $\Omega_{max} = Ma^2$ occurs at $\mu_{max} = \pm a\sqrt{\frac{2M}{\omega''}}$.

It is useful to write explicitly the form of the perturbed, small amplitude wave profile in the
 original coordinates $(x,t)$. Making the proper substitutions including an explicit addition of a small parameter $\epsilon$ to coincide with the derivation in the previous section, we obtain
\begin{align}
\eta(x,t) = \text{Re}\left[a e^{iM a^2 \epsilon^2 t}e^{i(kx-\omega t)} + u\delta e^{iM a^2 \epsilon^2 t}e^{i(kx-\omega t)} e^{\Omega \epsilon^2 t}e^{i\mu\epsilon(x-v_g t)}\right].
\label{eq:modPert}
\end{align}
We can interpret the unperturbed wave profile (wave profile of O($\delta^0$)) as a cosine solution of wavenumber $k$ that travels with speed $c_{\text{NLS}} = \omega - Ma^2\epsilon^2$. This gives the form of the time dependence of the perturbation as
\begin{align}
\bar{\Omega} - i\bar{\mu}v_g - ic_{\text{NLS}},
\end{align}
and the factor multiplying the spatial dependence in the exponential is $ik + i\bar{\mu}$, where the barred variables contain their respective powers of $\epsilon$. The bar notation will now be dropped for ease.


\section{Numerical Setup} \label{sec:numerics}

In this section, we describe how the solutions to  the Euler's equations given by \eqref{eq:euler} are computed and proceed by setting up the eigenvalue problem used to compute their stability. To do this, we first  switch into a travelling frame of reference, moving at speed $c$. This introduces a natural parametrisation of the problem with respect to the wave speed. We start by obtaining the solution to the linearised equations. Then using a continuation method, as we change the wave speed, we will compute more and more nonlinear solutions to the equations in the travelling frame. For each of these solutions, we examine their stability in a spectral sense as defined in the second part of this section. For more details on how these computations are done, see \cite{DO11, DT14}. The numerical results will be shown in subsequent sections.

\subsection{Numerical Solutions} \label{subsec:numSolutions}
We use the reformulation due to AFM for our numerical procedure. First, we rewrite the equations \eqref{eq:local} and \eqref{eq:nonlocal} in a travelling frame of reference moving at a constant speed $c$ with $x \rightarrow x-ct$
\begin{align}
 q_t - cq_x + \frac{1}{2} q_x^2 + g\eta -
\frac{1}{2}\frac{(\eta_t-c\eta_x+\eta_xq_x)^2}{1+\eta_x^2} = - D
P_{\text{flex}}
\label{eq:localC} \\
 \int_0^Le^{ikx}\left( i(\eta_t - c\eta_x) \cosh(k(\eta+h)) +
q_x\sinh(k(\eta+h)) \right)dx = 0 ~~ \ k \in \mathbb{Z}.
\label{eq:nonlocalC}
\end{align}
We now look for solutions that are stationary in the travelling frame of reference.  From the local equation \cite{DO11}, we obtain 
\begin{align}
 q_x = c - \sqrt{(1+\eta_x^2)\left(c^2 - 2g\eta - 2
DP_{\text{flex}}\right)}.
\label{eq:qx}
\end{align}
Using the form of $q_x$ in the nonlocal equation, we obtain one equation for the unknown wave profile $\eta(x)$, parametrised by the wave speed $c$
\begin{align}
  \int_0^{2\pi}e^{ikx} \sqrt{(1+\eta_x^2)\left(c^2 - 2g\eta - 2
DP_{\text{flex}}\right)}\sinh(k(\eta+h))dx = 0.
\label{eq:steady}
\end{align}
Alternatively,
\begin{align}
\int_0^{{2\pi}} e^{ikx}  \sqrt{(1+\eta_x^2) \left(c^2 - 2g\eta - 2D P_{\text{flex}}\right)} \left(\sinh(k\eta) + \cosh(k\eta)\tanh(kh)\right) dx = 0,~~ \ k \in \mathbb{Z},
\label{eq:tanhSS}
\end{align}
where we have separated the explicit dependence on the depth $h$. In the limit as $h\rightarrow \infty$, this gives
\begin{align}
\int_0^{{2\pi}} e^{ikx}  \sqrt{(1+\eta_x^2) \left(c^2 - 2g\eta - 2 D P_{\text{flex}}\right)} e^{|k|\eta} dx = 0,~~ \ k \in \mathbb{Z}.
\label{eq:infSS}
\end{align}
We can show that the solution is be symmetric for small amplitude \cite{TDW16} and therefore has the following perturbation series expansion
\begin{align}\label{eq:stokes}
  \eta(x) = 2\epsilon \cos x+ \sum_{k=2}^{\infty} \epsilon^k \eta_k(x), \quad
  \eta_k(x)=\sum_{j=2}^k 2\hat \eta_{kj}\cos(jx).
\end{align}
This is similar to the expansion introduced in Section \ref{sec:asymptotics}, equation \eqref{eq:waveDecomp}, but with the dependence on the small parameter $\epsilon$ shown explicitly.  Using this expansion, we notice that for $g$ and $h$ fixed, if the  flexural rigidity parameter $D$ satisfies
\begin{align}
  (g+D) K \tanh(h)  - \left( g  + K^4 D \right)\tanh(Kh) = 0, \qquad
  (K\neq 1),
\label{eq:resCond}
\end{align}
known as the resonance condition, $\eta_{K}$ will have vanishing denominators. This is similar to what happens for gravity-capillary waves \cite{VBbook} and leads to large values for coefficients of certain modes. However, these denominators never fully vanish due to the presence of nonlinearity when solutions are computed numerically.

Equation \eqref{eq:stokes} implies that at linear order, the solution is made up of one cosine mode of amplitude $\epsilon$. Numerically, we compute solutions in a cosine basis with the small parameter absorbed into the coefficients of the modes. Introducing a truncated series expansion for the wave profile as
\begin{align}
\eta_N(x) = \sum_{j=1}^N a_j\cos(jx)
\label{eq:numStokes}
\end{align}
with the number of modes given by $N$. Equation \eqref{eq:tanhSS} is valid for every integer $k\neq 0$. We let $k$ to take values from 1 to $N$. The $N$ equations we obtain have a free parameter $c$.  In practice, since the largest coefficient of the expansion is $a_1(=2\epsilon)$,  it  is the coefficient we use to control the wave amplitude along the bifurcation branch. This implies that for each point on the bifurcation branch, the vector of unknowns is given by
\begin{align*}
z = [c, a_2, a_3, \hdots, a_N]^T.
\end{align*}
The $N$ equations are then given by
\begin{align}
F^{(N)}_m(z) &= \int_0^{2\pi} e^{imx}\sqrt{\left(1+\eta_{N,x}^2\right) \left(c^2 - 2g\eta_N- 2 D P_{\text{flex}} \right)}(\sinh(m\eta_N)+\cosh(m\eta_N)\tanh(mh))dx,
\label{eq:numEq}
\end{align}
with $m = 1..N$. We wish to solve $F^{(N)}(z)=0$ for the unknown vector $z$. Using Newton's method, the $n$-th iteration is given by
\begin{align*}
z^{n+1} = z^{n} - J^{-1}(z^n)F(z^n),
\end{align*}
with $J$ the Jacobian. We now start the continuation method by noting that flat water ($\eta(x) = 0$) can travel at any wave speed $c$ as shown by \eqref{eq:tanhSS}. For one particular value of wave speed, we obtain a nontrivial solution $\eta(x) = a_1 \cos(x)$ with $a_1 = 2\epsilon$ small. This bifurcation point given by
\begin{align*}
[ \sqrt{\tanh(h)(g+D)}, 0, 0, 0, \hdots, 0]^T,
\end{align*}
is used as a start for a branch of solutions with increasing amplitudes. We substitute this guess into the equations given by \eqref{eq:numEq} to compute the actual profile of the wave using Newton's method. Then we scale up the wave amplitude and use this as a guess in a new step of Newton's method to compute a larger wave. Matlab is used to implement the numerical scheme. To check the convergence of the algorithm, we check the decay in the Fourier modes of the solution. If the decay is not sufficient, the number of Fourier modes is increased ensuring that the computed modes with highest wave number have negligible amplitude. This can be done for regular or resonant solutions and we will show these results in Section \ref{sec:genResults}.

\subsection{Numerical Stability} \label{subsec:numStability}
Ultimately, we are interested in analysing how the computed solutions behave if they are perturbed by a time-dependent perturbation. In this section, we follow the methodology for Fourier-Floquet-Hill (FFH) method outlined in \cite{DK06, DO11, DT14}.  So far, we have computed the travelling wave solutions which we will refer to as $\eta^0$ and we use \eqref{eq:qx} to compute the corresponding equilibrium derivative of the velocity potential $q^0_x$. We can now introduce the following perturbation with a particular time dependence as
\begin{align}
q({x},t) = q^{(0)}(x) + \delta q^{(1)}({x})e^{{\lambda} t} + O(\delta^2),\nonumber\\
\eta({x},t) = \eta^{(0)}({x}) + \delta \eta^{(1)}({x})e^{{\lambda} t} + O(\delta^2)
\label{eq:perts1}
\end{align}
where we are still in the travelling frame of reference. If $\lambda$ has some real, positive part, then the solution is exponentially growing and therefore it is not spectrally stable. 

Equations \eqref{eq:localC} and \eqref{eq:nonlocalC} are time dependent equations, with the second only valid for solutions of period $2 \pi$.  We do not want to restrict the period of the perturbations $q^{(1)}$ and $\eta^{(1)}$, which is possible by using Floquet's Theorem \cite{InceBook, DK06}. For our problem, this implies that perturbations   bounded in space may be decomposed as
\begin{align} 
q^{(1)}(x)  = e^{i\mu x}\tilde{q}^{(1)},~~
\eta^{(1)}(x)  = e^{i\mu x}\tilde{\eta}^{(1)},
\label{eq:perts2}
\end{align}
where $\mu\in [-1/2,1/2)$ is the Floquet exponent and $\tilde q^{(1)}$, $\tilde \eta^{(1)}$ are periodic with period $2\pi$. It is straightforward to apply the Floquet Theorem to the local equation, but the nonlocal case requires modification. We need to replace the integral in the nonlocal equation over one period by the average over the whole line
\begin{align}
\left<f(x)\right> = \lim_{M\rightarrow\infty}\frac{1}{M}\int_{-M/2}^{M/2}f(x)dx,
\end{align}
which is defined for almost periodic $f(x)$. This  includes (quasi-) periodic $f(x)$ as in (\ref{eq:perts2}) \cite{BohrBook}. Then we linearise the following system of equations about a traveling wave solution
\begin{align}
q_t - cq_x  + \frac{1}{2}q_x^2 + g\eta - \frac{1}{2}\frac{(\eta_t-c\eta_x+q_x\eta_x)^2}{1+\eta_x^2} & = D P_{\text{flex}}, \\
\lim_{M\rightarrow \infty} \frac{1}{M} \int_{-M/2}^{M/2} e^{ikx}\left[i(\eta_t-c\eta_x)\cosh(k(\eta+h))+q_x\sinh(k(\eta+h)) \right] dx & = 0, \ \ k\in \Lambda.
\end{align}
Using \eqref{eq:perts1}, ignoring terms of $O(\delta^2)$ and higher, and dropping the tildes, we obtain
\begin{align}
\lambda\left(f\eta^{(1)}-q^{(1)}\right) & = (q^{(0)}_{x}-c)\mathrm{D}_xq^{(1)} + g\eta^{(1)}-f\left[ (q^{(0)}_{x}-c)\mathrm{D}_x\eta^{(1)} + \eta^{(0)}_{x}\mathrm{D}_xq^{(1)}\right] + f^2\eta^{(0)}_{x}\mathrm{D}_x\eta^{(1)}\nonumber + DG(\eta^{(0)},\eta^{(1)})  \\
\lambda\left<e^{ikx}\left[-i\mathcal{C}_k\eta^{(1)}\right]\right>&= \left<e^{ikx}\left[-i\mathcal{C}_kc\mathrm{D}_x\eta^{(1)} + \mathcal{S}_k\mathrm{D}_xq^{(1)}+\left(-i\eta^{(0)}_{x}c\mathcal{S}_k+ q^{(0)}_{x}\mathcal{C}_k\right)k\eta^{(1)}\right]\right>,
\label{eq:linNL}
\end{align}
where
\begin{align*}
f(\eta^{(0)},q^{(0)}) &= \frac{\eta^{(0)}_{x}(q^{(0)}_{x}-c)}{1+(\eta^{(0)}_{x})^2},& \mathrm{D}_x & =i\mu+\partial_x, & \\
\mathcal{S}_k &= \sinh(k(\eta^{(0)} + h)), & \mathcal{C}_k & = \cosh(k(\eta^{(0)}+h)), & \mathcal{T}_k = \tanh(k(\eta^{(0)}+h)).
\end{align*}
The term $G(\eta_1, q_1)$ depends on the model we use for the waves under ice. For example, for the nonlinear (Toland) model, this is given by
\begin{align*}
G(\eta^{(1)}, q^{(1)}) & = \frac{\eta^{(1)}_{4x}}{(1+(\eta^{(0)}_{x})^2)^{5/2}} - \frac{5\eta^{(0)}_{4x}\eta^{(0)}_{x}\eta^{(1)}_{x}}{(1+(\eta^{(0)}_{x})^2)^{7/2}}
-10\frac{\eta^{(1)}_{x}\eta^{(0)}_{2x}\eta^{(0)}_{3x}}{(1+(\eta^{(0)}_{x})^2)^{7/2}} - 10\frac{\eta^{(0)}_{x}\eta^{(1)}_{2x}\eta^{(0)}_{3x}}{(1+(\eta^{(0)}_{x})^2)^{7/2}} \\
& -10\frac{\eta^{(0)}_{x}\eta^{(0)}_{2x}\eta^{(1)}_{3x}}{(1+(\eta^{(0)}_{x})^2)^{7/2}} + 70 \frac{(\eta^{(0)}_{x})^2\eta^{(0)}_{2x}\eta^{(0)}_{3x}\eta^{(1)}_{x}}{(1+(\eta^{(0)}_{x})^2)^{9/2}} - \frac{15}{2}\frac{(\eta^{(0)}_{2x})^2\eta^{(1)}_{2x}}{(1+(\eta^{(0)}_{x})^2)^{9/2}} + \frac{45}{2}\frac{(\eta^{(0)}_{2x})^3\eta^{(0)}_{x}\eta^{(1)}_{x}}{(1+(\eta^{(0)}_{x})^2)^{11/2}} \\
& + 30\frac{(\eta^{(0)}_{2x})^3\eta^{(0)}_{x}\eta^{(1)}_{x}}{(1+(\eta^{(0)}_{x})^2)^{9/2}}+45\frac{(\eta^{(0)}_{x})^2(\eta^{(0)}_{2x})^2\eta^{(1)}_{2x}}{(1+(\eta^{(0)}_{x})^2)^
{9/2}}-135\frac{(\eta^{(0)}_{x})^3(\eta^{(0)}_{2x})^3\eta^{(1)}_{x}}{(1+(\eta^{(0)}_{x})^2)^{11/2}}.
\end{align*}
\noindent Since $q^{(1)}$ and $\eta^{(1)}$ are periodic with period $2\pi$,
\begin{align}
{q}^{(1)} = \sum_{m=-\infty}^{\infty} {Q}_m e^{imx}, ~~
{\eta}^{(1)} = \sum_{m=-\infty}^{\infty} {N}_m e^{imx},
\end{align}
\noindent with

\begin{align}
{Q}_n =  \frac{1}{2\pi}\int_0^{2\pi} e^{-inx} {q}^{(1)}(x) dx,~~~~
{N}_n =  \frac{1}{2\pi}\int_0^{2\pi} e^{-inx} {\eta}^{(1)}(x) dx.
\end{align}

Truncating to the $2N+1$ Fourier modes from $-N$ to $N$, we define the unknowns as

\begin{align}
U(x) = \begin{bmatrix} {N}_{-N}(x), \hdots, {N}_0(x),  \hdots, {N}_N(x), {Q}_{-N}(x), \hdots {Q}_0(x), \hdots, {Q}_N(x) \end{bmatrix}^T.
\end{align}

\noindent This leads to the finite-dimensional generalised eigenvalue problem

\begin{align}
\lambda \mathcal{L}_1 U(x) = \mathcal{L}_2 U(x)
\label{eq:matrix}
\end{align}
\noindent where

\begin{align}
\mathcal{L}_1 = \begin{bmatrix} {A} & {-I} \\ {C} & {0} \end{bmatrix}, ~~\mathcal{L}_2 = \begin{bmatrix} {S} & {T} \\ {U} & {V} \end{bmatrix}
\label{eq:BM}
\end{align}

\noindent with $I$ and $0$ the $(2N+1)\times(2N+1)$ identity and zero matrix, respectively. The blocks $A$, $S$ and $T$ originate from the local equation, while $C, U$ and $V$ come from the nonlocal equation. The matrix entries are given by

\begin{align*}
{A}_{m,n} &=  \frac{1}{2\pi} \int_0^{2\pi} e^{i(m-n)x} f dx, ~~
{C}_{m,n} = -i  \frac{1}{2\pi} \int_0^{2\pi} e^{i(m-n)x}\mathcal{C}_{\mu+m} dx,\\
{S}_{m,n} & = -  \frac{1}{2\pi} \int_0^{2\pi} e^{i(m-n)x} \left[ - g + f(q^{(0)}_{x} -c)i(\mu+(m-N)) - f^2\eta^{(0)}_{x} i(\mu+(m-N)) + D G(\cdot) \right] dx, \\
{T}_{m,n} & =  \frac{1}{2\pi} \int_0^{2\pi} e^{i(m-n)x} \left[ (q^{(0)}_{x}-c)i(\mu+(m-N)) - f \eta^{(0)}_{x}i(\mu+(m-N)) \right] dx, \\
{U}_{m,n} & =  \frac{1}{2\pi} \int_0^{2\pi} e^{i(m-n)x} \left[S_{\mu+m} i(\mu+(m-N)) \right] dx, \\
{V}_{m,n} & =  \frac{1}{2\pi} \int_0^{2\pi} e^{i(m-n)x} \left[-ic(\mu+(m-N)) \mathcal{C}_k + k (-i \eta^{(0)}_{x} c \mathcal{S}_{\mu+m}) + q^{(0)}_{x} \mathcal{C}_{{\mu+m}} \right] dx.
\end{align*}

\noindent Lastly,

\begin{align*}
\mathcal{C}_{{\mu+m}} = \cosh(({\mu+m})\eta^{(0)}) + \mathcal{T}_{{\mu+m}}\sinh(({\mu+m})\eta^{(0)}), ~~ \mathcal{S}_{{\mu+m}} = \sinh(({\mu+m})\eta^{(0)}) + \mathcal{T}_{{\mu+m}}\cosh(({\mu+m})\eta^{(0)}),
\end{align*}

\noindent with $\mathcal{T}_{{\mu+m}} = \tanh(({\mu+m})h)$. All block matrices in (\ref{eq:BM}) are of size $(2N+1)\times(2N+1)$ with $N$ the number of modes we retain. The convergence properties of the Floquet-Fourier-Hill method (FFH) as $N\rightarrow \infty$ are discussed in \cite{CD10, JZ12}. 

In order to compare with the result of the previous section, we need to compare the unstable perturbation in a stationary frame of reference. If we substitute the transformations into \eqref{eq:perts1}, we obtain that the perturbed surface elevation is now given by
\begin{align}
\eta(x,t) = \eta_0(x-ct) + \delta \text{Re}\left[e^{\lambda t} \sum_{m=-N}^{N} \hat{N}_m e^{i(m+\mu)(x-ct)}\right]
\label{eq:hillPert}
\end{align}
where we note that $\hat{N}_m$ may be complex.

\section{Results in the NLS Regime}\label{sec:NLSresults}
In this section, we show how the asymptotic results and numerical results coincide in the same regime. We start by  examining solutions to  the Euler's equations and then by discussing their stability. We do this for water of infinite depth. We focus on 5 different regimes summarised in Table \ref{tab:summary}. 
\begin{table}[h]
\begin{center}
\begin{tabular}{ | c | c | c | c | c |}
\hline
Regime &  \multicolumn{2}{|c|}{Bifurcation Branch Direction} & \multicolumn{2}{|c|}{Modulational Instability} \\
\hline
Deep Water ($h=\infty$) & Linear & Nonlinear & Linear & Nonlinear \\ \hline
$D = 0.01$ & \textcolor{blue}{right} & \textcolor{red}{right} & \textbf{\textcolor{blue}{unstable}} & \textbf{\textcolor{red}{unstable}} \\
$D = 0.05$ & \textcolor{blue}{right} & \textcolor{red}{right} & \textcolor{blue}{stable} & \textcolor{red}{stable}\\
$D = 0.1$ & \textcolor{blue}{left} & \textcolor{red}{left} & \textbf{\textcolor{blue}{unstable}} & \textbf{\textcolor{red}{unstable}}\\
$D = 0.3$ & \textcolor{blue}{right} & \textcolor{red}{right} & \textcolor{blue}{stable} & \textcolor{red}{stable}\\
\hline
$D = 25$ & \textcolor{blue}{right} & \textcolor{red}{left} & \textbf{\textcolor{blue}{unstable}} & \textcolor{red}{stable}\\
\hline
\end{tabular}
\end{center}
  \caption{Summary of the results for deep water showing the direction of the bifurcation branch of solutions and which solutions exhibit modulational instabilities for both linear and nonlinear models in different flexural rigidity regimes.\label{tab:summary}}
\end{table}

\subsection{Solutions} \label{subsec:NLSsolutions}

\begin{figure}[h!]
\begin{center}
\begin{tikzpicture}
    \node[anchor=south west,inner sep=0] at (0,0) {\includegraphics[width=0.68\textwidth]{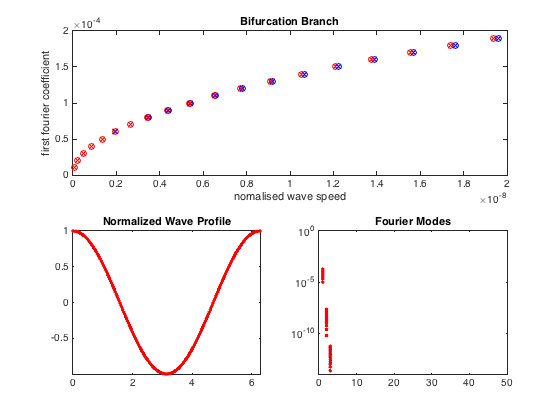}};
    \node at (4.5, 6.03) {LIN};
    \node at (4.2, 6.63) {NL};
\end{tikzpicture}
\end{center}
\caption{Solutions with $D=0.01$ and $h=\infty$. Top panel shows the bifurcation branch with circles the numerical computations and crosses showing the NLS approximation. The linear model (blue, labelled LIN) extending slightly further than the nonlinear (red, labelled NL). Bottom panel shows the profile (left) and semilog plot of the Fourier coefficients (right). Few Fourier modes imply we are close to the bifurcation point and the profiles look the same for both models.
\label{fig:deepSolnsD0_01}}
\end{figure}

\begin{figure}[h!]
\begin{center}
\begin{tikzpicture}
    \node[anchor=south west,inner sep=0] at (0,0) {\includegraphics[width=0.68\textwidth]{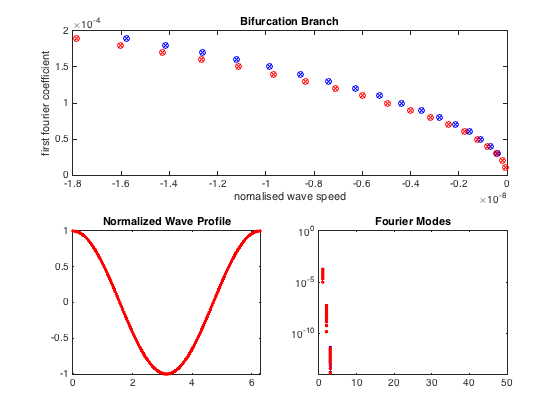}};
    \node at (4.4, 7.2) {LIN};
    \node at (4.2, 6.63) {NL};
\end{tikzpicture}
\end{center}
\caption{Solutions with $D=0.1$ and $h=\infty$. Top panel shows the bifurcation branch with circles the numerical computations and crosses the NLS approximation. The nonlinear (Toland) model for ice (red, labelled NL) extending further to the left than the linear (biharmonic) model (blue, labelled LIN). Botton panel shows the profile (left) and semilog plot of the Fourier coefficients (right). Few Fourier modes imply we are close to the bifurcation point and the profiles look the same. \label{fig:deepSolnsD0_1}}
\end{figure}

\begin{figure}[h!]
\begin{center}
\begin{tikzpicture}
    \node[anchor=south west,inner sep=0] at (0,0) {\includegraphics[width=0.68\textwidth]{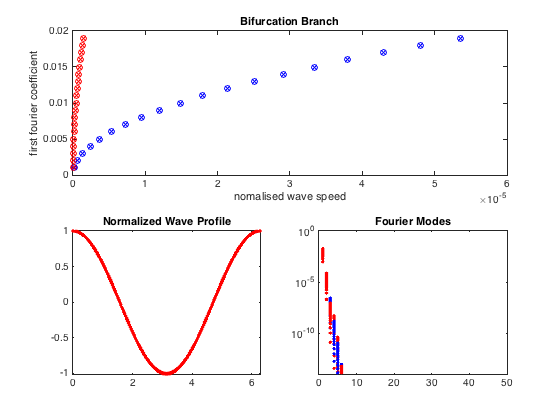}};
    \node at (1.95, 6.63) {NL};
    \node at (4.0, 6.63) {LIN};
\end{tikzpicture}
\end{center}
\caption{Solutions with $D=0.3$ and $h=\infty$. Top panel shows the bifurcation branch with circles the numerical computations and crosses the NLS approximation. Linear (biharmonic) model for ice extending further to the right (blue, labelled LIN) than the nonlinear (Toland) model (red, labelled NL). Botton panel shows the profile (left) and semilog plot of the Fourier coefficients (right). Few Fourier modes imply we are close to the bifurcation point and the physical profiles are same for both models.  \label{fig:deepSolnsD0_3}}
\end{figure}

\begin{figure}
\begin{center}
\begin{tikzpicture}
    \node[anchor=south west,inner sep=0] at (0,0) {\includegraphics[width=0.68\textwidth]{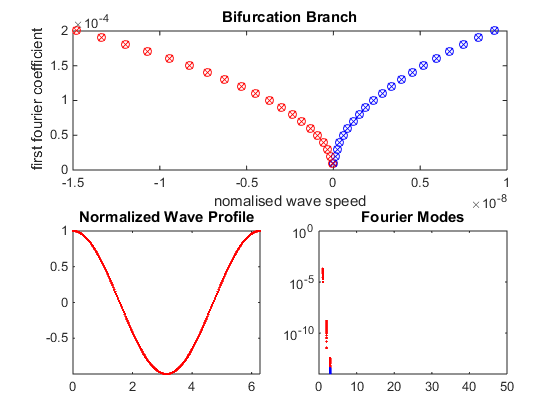}};
    \node at (5.32, 5.63) {NL};
    \node at (7.28, 5.63) {LIN};
\end{tikzpicture}
\end{center}
\caption{Solutions with $D=25$ and $h=\infty$. Top panel shows the bifurcation branch with circles the numerical computations and crosses the NLS approximation. Linear (biharmonic) model for ice extending to the right (blue, labelled LIN) and the nonlinear (Toland) model to the left (red, labelled NL). Botton panel shows the profile (left panel) and semilog plot of the Fourier coefficients (right panel). Few Fourier modes mean we are in the linear regime with two models giving the same physical profiles.\label{fig:deepSolnsD25}}
\end{figure}

The numerical results for  different values of flexural rigidity $D$ are shown in Figures \ref{fig:deepSolnsD0_01} - \ref{fig:deepSolnsD25}. We use the convention of the linear model of elasticity (biharmonic) in blue and the nonlinear model (Toland) in red also labelled as NL and LIN respectively. We computed these solutions for five distinct values of the flexural rigidity, focussing on the regions for which we have different stability results according to the NLS derivation as shown in Figure \ref{fig:NLSCoeff} and summarised in Table \ref{tab:summary}. To check how nonlinear these solutions are, we compare them to a bifurcation branch we get from the NLS approximation given by \eqref{eq:modPert}. As we have shown,
\begin{align}
c_{\text{NLS}} = \sqrt{1+D} - M{a^2}.
\label{eq:cNLS}
\end{align}
If we assume that the  waves are of period $2\pi$ as was done for  the numerical solutions, then  the solutions in the NLS regime will be well approximated  by a cosine with $k=1$ and amplitude $a/2$. We plot the amplitude $a$ and $c_{\text{NLS}}$ using crosses and the fully nonlinear results obtained from the procedure outlined in Section \ref{subsec:numSolutions} using circles for waves in  infinite depth ($h \rightarrow \infty)$. Since the goal is to compare numerical stability results to asymptotic results from NLS, the solutions for which we analyse stability should stay close to those given by \eqref{eq:cNLS} which are approximated by a cosine. These solutions are shown in Figures \ref{fig:deepSolnsD0_01}-\ref{fig:deepSolnsD25}. In these figures, the top panel shows the bifurcation branch with the normalised wave speed where we subtract the speed at the bifurcation and the normalised wave profile computed by dividing by the maximum amplitude of the wave shown in the bottom left and the semilog plot of the corresponding Fourier coefficients shown on the bottom right. We note that the computations use 50 coefficients, but only a few modes are needed  for low amplitude waves. These figures each show the bifurcation branch for which the numerical solutions and the asymptotic solutions overlap. While we can see a difference in wave speed for the two models for the  ice, we do not  see this in the normalised profiles shown in the lower panel on the left. The Fourier modes for these profiles shown in the right panel are very similar for both models as well.  Figure \ref{fig:deepSolnsD0_01} for solutions with $D=0.01$ and Figure \ref{fig:deepSolnsD0_1} for $D=0.1$ show that the bifurcation branch direction is different in these two regimes. However, both models and both regimes are well approximated by the NLS as they contain few Fourier modes.  As the values of the flexural rigidity is increased, we see that the models give different solutions as illustrated in Figures \ref{fig:deepSolnsD0_3} and \ref{fig:deepSolnsD25}, with the latter showing that depending on the model for flexural-gravity waves, the bifurcation branches change directions. This implies that in the linear model, high amplitude waves travel faster than lower amplitude whereas the nonlinear model is the opposite.

\subsection{Stability Results in the NLS Regime} \label{subsec:NLSstability}

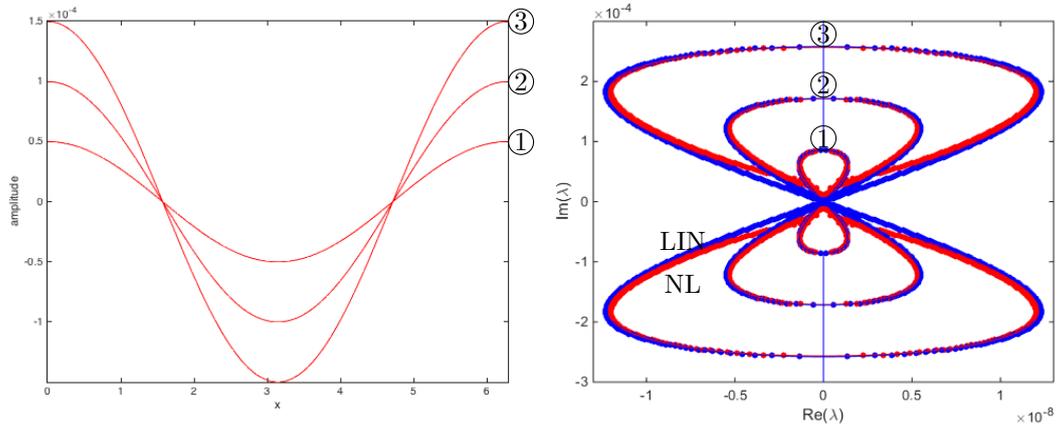
\begin{figure}[h!]
\begin{center}
\begin{tikzpicture}
    \node[anchor=south west,inner sep=0] at (0,0) {\includegraphics[width=0.49\textwidth]{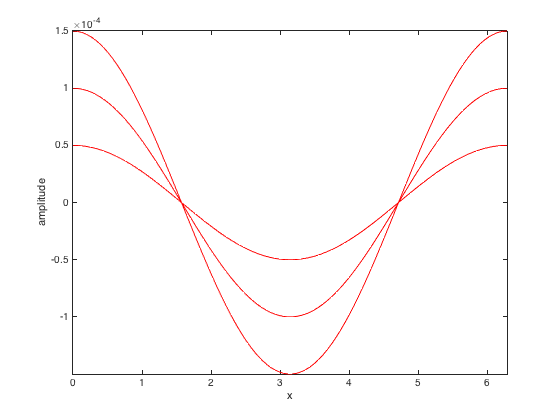}};
\node[anchor=south east, inner sep=0] at (15,0) 
{\includegraphics[width=0.49\textwidth]{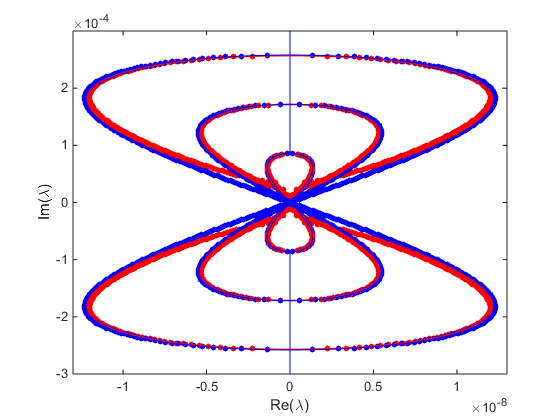}};
	\node at (9.38, 2.52) {LIN};
    \node at (9.38, 1.95) {NL};
    \draw (7.25,5.43) circle (0.165cm);
    \node at (7.25, 5.43) {3};
    \draw (7.25,4.63) circle (0.165cm);
    \node at (7.25, 4.63) {2};
    \draw (7.25,3.83) circle (0.165cm);
    \node at (7.25, 3.83) {1}; 
    \draw (11.23,3.88) circle (0.165cm);
    \node at (11.23,3.88) {1};
    \draw (11.23,4.59) circle (0.165cm);
    \node at (11.23,4.59) {2};
    \draw (11.23,5.26) circle (0.165cm);
    \node at (11.23,5.26) {3};
\end{tikzpicture}
\end{center}
\caption{The regime where $D=0.01$, infinitely deep water. On the left are the wave profiles for which we see the complex eigenvalue plane on the right. In blue is the linear model (labelled LIN) with curves lying outside the nonlinear model in red (labelled NL). \label{fig:stabD0_01}}
\end{figure}

\begin{figure}[h!]
\begin{center}
\begin{tikzpicture}
    \node[anchor=south west,inner sep=0] at (0,0) {\includegraphics[width=0.49\textwidth]{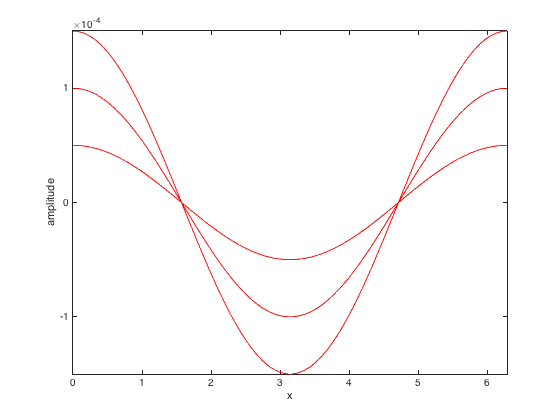}};
\node[anchor=south east, inner sep=0] at (15,0) 
{\includegraphics[width=0.49\textwidth]{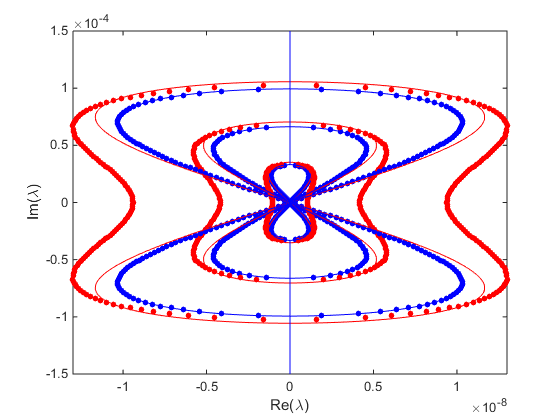}};
	\node at (9.3, 1.9) {LIN};
    \node at (9.15, 1.2) {NL};
    \draw (7.25,5.43) circle (0.165cm);
    \node at (7.25, 5.43) {3};
    \draw (7.25,4.63) circle (0.165cm);
    \node at (7.25, 4.63) {2};
    \draw (7.25,3.83) circle (0.165cm);
    \node at (7.25, 3.83) {1}; 
    \draw (11.23,3.58) circle (0.165cm);
    \node at (11.23,3.58) {1};
    \draw (11.23,4.10) circle (0.165cm);
    \node at (11.23,4.10) {2};
    \draw (11.23,4.62) circle (0.165cm);
    \node at (11.23,4.62) {3};
\end{tikzpicture}
\end{center}
\caption{The regime where $D=0.1$ in infinitely deep water. On the left are the wave profiles for which we see the complex eigenvalue plane on the right. In blue is the linear (labelled LIN) model with curves lying inside the nonlinear (NL) model in red.\label{fig:stabD0_1}}
\end{figure}
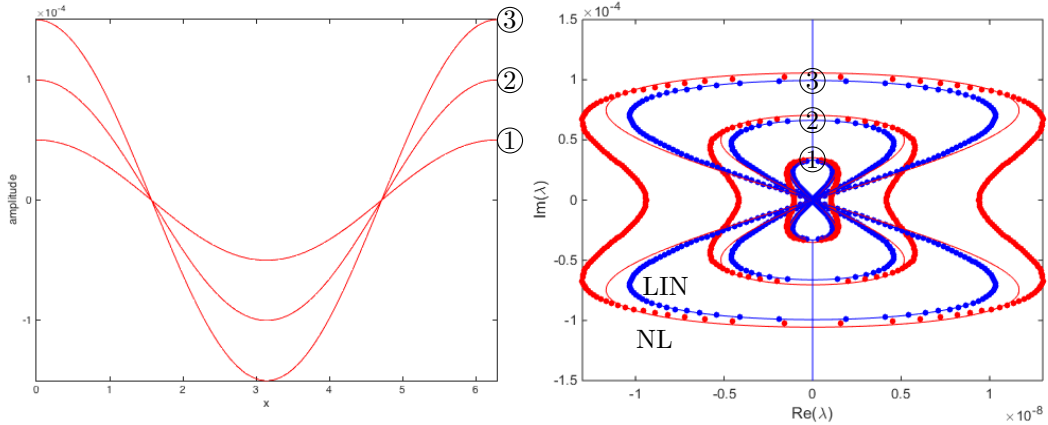

\begin{figure}
\begin{center}
\begin{tikzpicture}
\node[inner sep=0] at (0,0) {\includegraphics[width=0.49\textwidth]{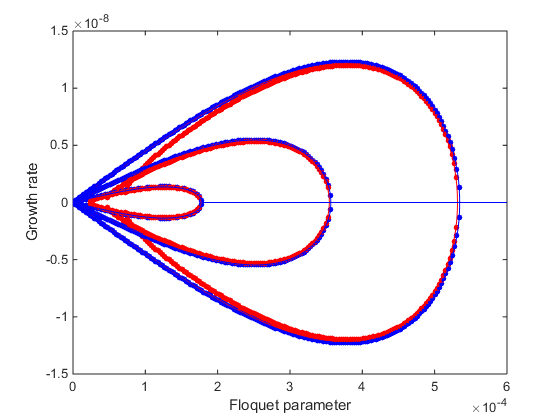}};
\node[inner sep=0] at (8,0) 
{\includegraphics[width=0.49\textwidth]{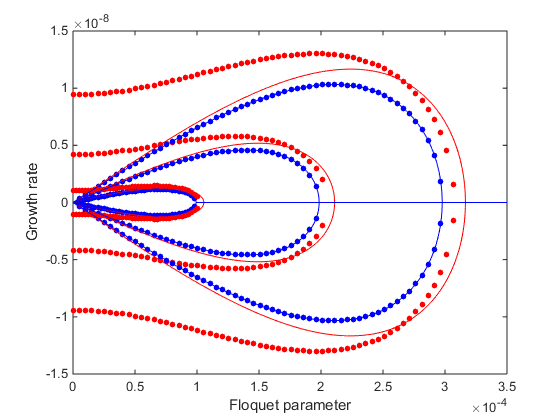}};
	\node at (0.15, 1.6) {NL};
    \node at (-0.3, 2.0) {LIN};
    \node at (10.65, 1.4) {NL};
    \node at (9.15, 1.4) {LIN};
    \draw (-0.88,0.1) circle (0.165cm);
    \node at (-0.88,0.1) {1};
    \draw (0.95,0.1) circle (0.165cm);
    \node at (0.95,0.1) {2};
    \draw (2.7,0.1) circle (0.165cm);
    \node at (2.7,0.1) {3};
    \draw (7.0,0.1) circle (0.165cm);
    \node at (7.0,0.1) {1};
    \draw (8.6,0.1) circle (0.165cm);
    \node at (8.6,0.1) {2};
    \draw (10.45,0.1) circle (0.165cm);
    \node at (10.45,0.1) {3};
\end{tikzpicture}
\end{center}
\caption{The perturbations leading to the largest instabilities for $D=0.01$ on the left and $D=0.1$ on the right in the same regime as waves in Figures \ref{fig:stabD0_01} and \ref{fig:stabD0_1}. In blue is the linear model (labelled LIN) and in red is the nonlinear model (labelled NL). The solid lines are the predictions via NLS and dotted lines are numerics. \label{fig:unstabPert}}
\end{figure}
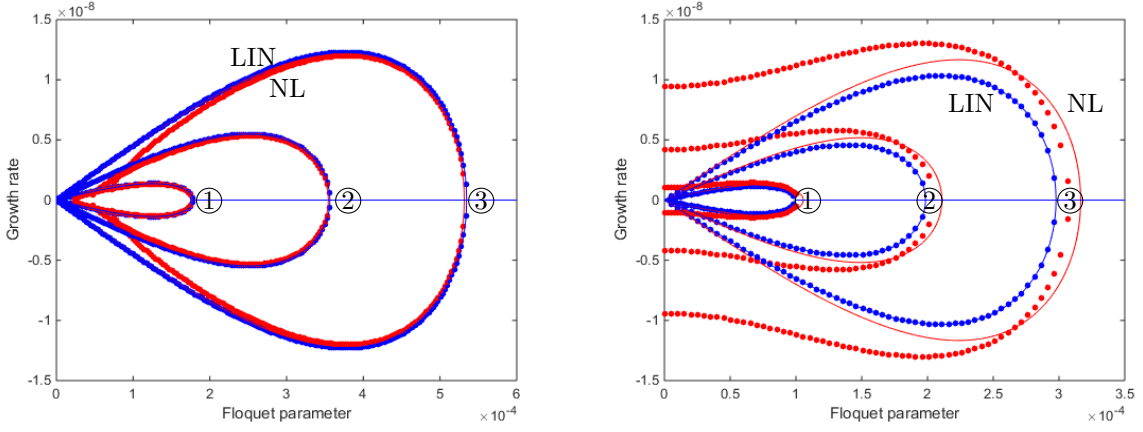

We proceed by analysing the stability of solutions computed above. We compare the modulational instability according to the asymptotic analysis through the NLS equation as seen in Section \ref{sec:asymptotics} with the numerical results from the method described  in Section \ref{subsec:numStability}, focussing on stationary waves of period $2\pi$, perturbed by a wave of any period. The asymptotic results assume that we are perturbing the mode $k=1$ with something that is of a similar wavenumber.  This implies that to compare, we need to set $m=\pm 1$ in \eqref{eq:hillPert}. For the full solution to be real, the resulting perturbed wave profile is of the form
\begin{align}
\eta_{\text{FFH}}(x,t) = \eta_0(x-ct) + \hat{N}_{1}e^{\lambda t} e^{-i\mu ct}e^{i(x-ct)}e^{i\mu x},
\label{eq:FFHModPert}
\end{align}
which we compare to the perturbation from the asymptotic method given by
\begin{align}
\eta_{\text{NLS}}(x,t) = a e^{iM a^2 \epsilon^2 t}e^{i(kx-\omega t)}  + \delta  u e^{\Omega t} e^{-i\mu v_g t}e^{i(x-c_{\text{NLS}}t)}e^{i\mu x}.
\label{eq:NLSModPert}
\end{align}
We are interested in how $\lambda$ and $\Omega$ compare. In examining the 5 regions outlined in Table \ref{tab:summary} numerically, we see stability where we anticipated, but we further examine the unstable regions for $D=0.01$ and $D=0.1$. We compute the stability spectrum of three different solutions and see that modulational instabilities are present for both models, as shown in Figures \ref{fig:stabD0_01} and \ref{fig:stabD0_1}. In these figures, the three solutions for which we analyse the instabilites are shown on the left and labelled 1 through 3, with solutions to both models overlapping. We see these resemble a cosine of different amplitudes. Their spectra is plotted on the right, with the corresponding labels. In these figures we plot $\text{Re}(\lambda)$ versus $\text{Im}(\lambda)$ as a series of points for all values of $\mu$ and for comparison, the asymptotic results are plotted as solid lines with $\Omega$ on the horizontal axis given by equation \eqref{eq:NLSGrowthRate} and on the vertical, $\mu(v_g -\omega)$. We see that for flexural-gravity waves modelled via the linear (biharmonic) model, the asymptotics and the numerics line up very well as shown in blue, but with the nonlinear (Toland) model, these deviate more, with the modulational instability beginning to change and move away from the origin in the spectral plane, as shown in red and labelled NL.  We also compare which perturbations lead to more unstable growth rates in Figures \ref{fig:unstabPert} for $D=0.01$ and $D=0.1$ with the smallest solutions giving the inner most set of results and the largest amplitude solutions giving the two outermost lines. We see that the lower coefficient of flexural rigidity, the nonlinear model gives a smaller result for the growth rate and for the larger coefficient, the model shows larger growth rates. Once again, the numerical and asymptotic results agree. 

\section{More General Results}\label{sec:genResults}
\begin{figure}[h]
\begin{center}
\includegraphics[width=0.49\textwidth]{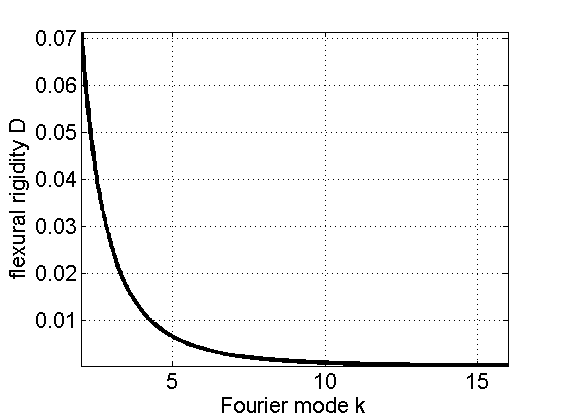}
\includegraphics[width=0.49\textwidth]{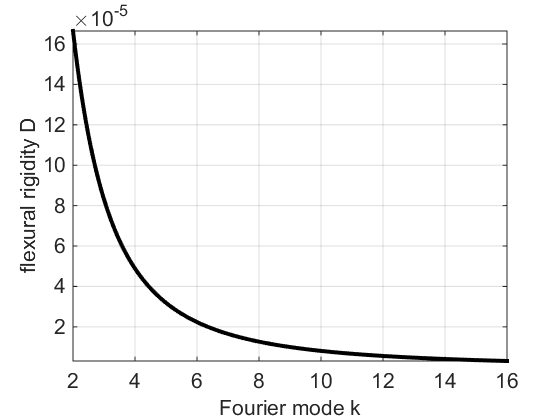}
\caption{Flexural rigidity $D$ as a function of wave number for which the resonance condition \eqref{eq:resCond} is satisfied. The left plot is for waves in infinite depth and the right is for $h=0.05$. \label{fig:resonance}}
\end{center}
\end{figure}

We examine the resonant regime. For the asymptotic regime governed by NLS with the nonlinear coefficient given by \eqref{eq:MTol} - \eqref{eq:MLin}, we see that the denominator blows up for $g-14k^4D=0$ or $D=1/14\approx0.07$ for $g=1$. This is a manifestation of resonance, which has been analysed for capillary-gravity waves and referred to as Wilton ripples \cite{W15, VBbook} and more recently by \cite{TDW16}. Outside of this regime, this condition is more generally given by \eqref{eq:resCond} and shown in Figure \ref{fig:resonance} for infinite depth and finite depth ($h=0.05$). We see that on the left of the figure, $D=0.07$ for $K=2$. These figures show that if we treat the flexural rigidity as a parameter, there will be a particular Fourier mode for which the resonance condition will hold, resulting in a large coefficient for that Fourier mode. We also note that the line $K=1$ is a vertical asymptote, which implies that the larger the coefficient of rigidity, the closer we get to the first mode being resonant.

\subsection{Resonant Solutions} \label{subsec:resonance}
\begin{figure}[h!]
\begin{center}
\includegraphics[width=0.68\textwidth]{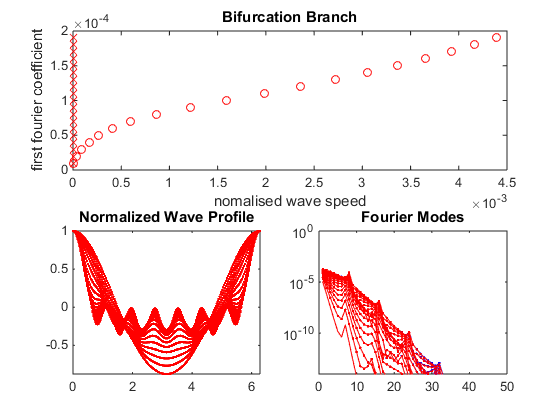}
\caption{Resonant solutions with $K=7$ and $h=0.05$. Top panel shows the bifurcation branch with circles the numerical computations and crosses the NLS approximation. Botton panel shows the profile (left) and semilog plot of the Fourier coefficients (right) showing resonance and the harmonics. Both models give the same result.\label{fig:K7Res}}
\end{center}
\end{figure}

\begin{figure}[h!]
\begin{center}
\includegraphics[width=0.68\textwidth]{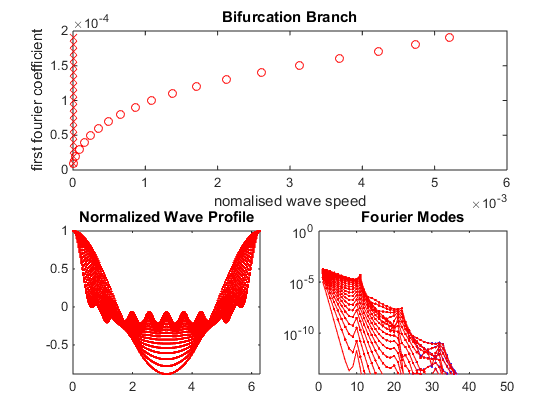}
\caption{Resonant solutions with $K=10$ and $h=0.05$. Top panel shows the bifurcation branch with circles the numerical computations and crosses the NLS approximation. Botton panel shows the profile (left) and semilog plot of the Fourier coefficients (right) showing resonance and the harmonics. Both models give the same result. \label{fig:K10Res}}
\end{center}
\end{figure}

In this section, we analyse what happens for solutions in the resonant regime in  water of finite
depth, using \eqref{eq:tanhSS}. We can rearrange the
formulation  in such a way that if we want resonance to occur at a particular wavenumber $k=K$, then we can set $D$ to satisfy \eqref{eq:resCond}. Figure \ref{fig:resonance} shows the flexural rigidity as a function of wavenumber $k$ for a nondimensional wave depth $h=0.05$. This depth was picked for illustrative purposes only and it can be compared with the results in \cite{TDW16} for capillary-gravity waves. For illustrative purposes, we pick the flexural rigidity parameter so that the resonant mode is $K = 7$ (i.e. $D \approx 1.65 \times 10^{-5}$) as presented  in Figure \ref{fig:K7Res}, which shows 7 secondary minima and the resonant mode $K=10$ (i.e. $D \approx 8.11 \times 10^{-6}$) as shown in Figure \ref{fig:K10Res} where we see 10 secondary minima in the bottom left part of the plot of the normalised wave profile. As before for infinite depth, we once again plot the NLS approximations as crosses and the numerical solutions to the full problem as zeros. In this regime, we are outside of the validity of the NLS approximation. However, the results for the two different models for the  ice are the same. In this case, we also see that more Fourier modes are needed to fully represent the solutions and  that they no longer decay exponentially but instead show humps at the resonant modes as well as the harmonics of those modes, particular at large amplitudes.

\subsection{High Frequency Instabilities} \label{subsec:highFrequency}
\begin{figure}
\begin{center}
\includegraphics[width=0.48\textwidth]{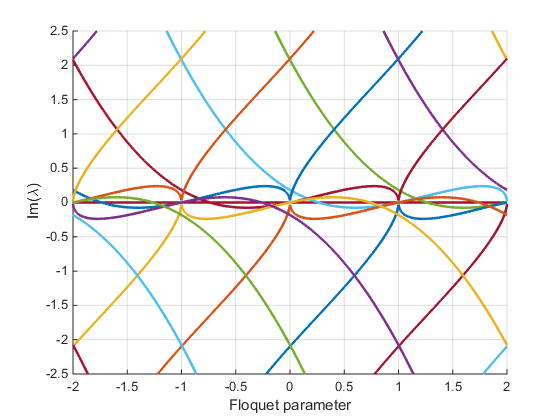}
\includegraphics[width=0.48\textwidth]{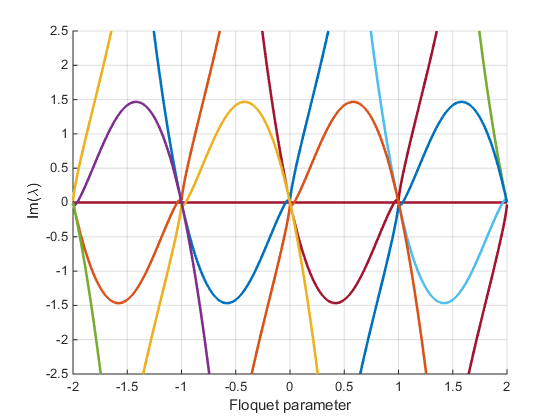}
\caption{Several eigenvalues given by \eqref{eq:evals} for $D=0.1$ on the left and $D=25$ on the right. For large values of $D$, there are more crossings at the origin. \label{fig:eCollisions}}
\end{center}
\end{figure}

Since the water wave problem is Hamiltonian \cite{Z68}, the spectra of any travelling wave solution is symmetric with respect to both the real and imaginary axes. Thus, in order for the solution to be spectrally stable, it is necessary for the spectrum to be on the imaginary axis, {\em i.e.}, $Re\{\lambda\} = 0$. It is well known that the eigenvalues corresponding to different Floquet exponents do not interact \cite{DK06}, thus we may restrict our attention to a fixed $\mu$ value.  These eigenvalues will depend on the solution to  the Euler's equations and in general, their analytic form is not known. However, we can compute them for a zero amplitude solution and they are given by
\begin{align}
\lambda _{\mu+m}^{\pm}= i c(\mu+m) \pm i \sqrt{\left[g(\mu+m) + D(\mu+m)^5\right]\tanh{((\mu+m)h)}}.
\label{eq:evals}
\end{align}
It is easy to see that these eigenvalues are on the imaginary axis and the flat water state is spectrally stable. The spectrum of (\ref{eq:linNL}) is a continuous function of the parameters appearing in $\mathcal{L}_1$ and $\mathcal{L}_2$ \cite{HislopAndSigalBook}, mainly the amplitude of the solution. In order for eigenvalues to leave the imaginary axis, they do so in pairs via eigenvalue collisions, which are a necessary condition for the development of instabilities \cite{MS86}. Thus we examine for which parameter values different eigenvalues shown in (\ref{eq:evals}) collide,
\begin{align}
\lambda^{s_1}_{\mu} = \lambda_{\mu+m}^{s_2} \ \ \text{for any} \ m\in\mathbb{Z}, s_1\neq s_2,
\label{eq:collision}
\end{align}
with $s_1$ and $s_2$ either positive or negative signs. We plot these eigenvalues for a particular set of parameters in the resonant regime. For this purpose, we unfold the Floquet parameter values to be outside of the usual range from $-0.5$ to $0.5$, effectively plotting several periods of the eigenvalues. Setting $D=0.1$, $h = \infty$, the eigenvalue collisions are shown on the left panel of Figure \ref{fig:eCollisions} and $D=25$ on the right. We see as $D$ is increased, more collisions are found closer to the origin, with a lot of eigenvalues meeting very close to the same value of Floquet parameter.
These collisions of eigenvalues may result in an instability that is different from a modulational instability. It is important to note that the resonance condition is equivalent to the collision condition for $\mu = 0$. This implies that resonant solutions should exhibit an instability near the origin of the complex eigenvalue plane.

\begin{figure}
\includegraphics[width=0.32\textwidth]{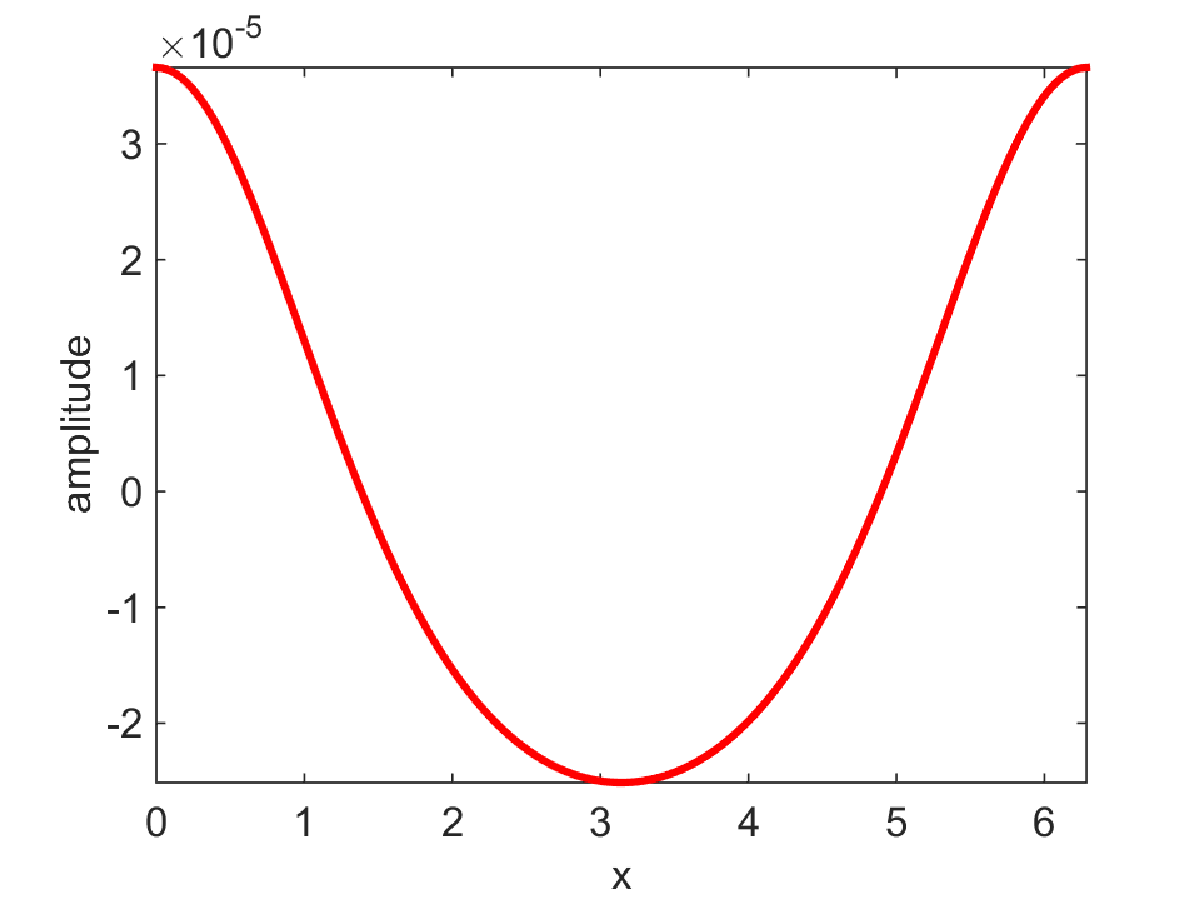}
\includegraphics[width=0.32\textwidth]{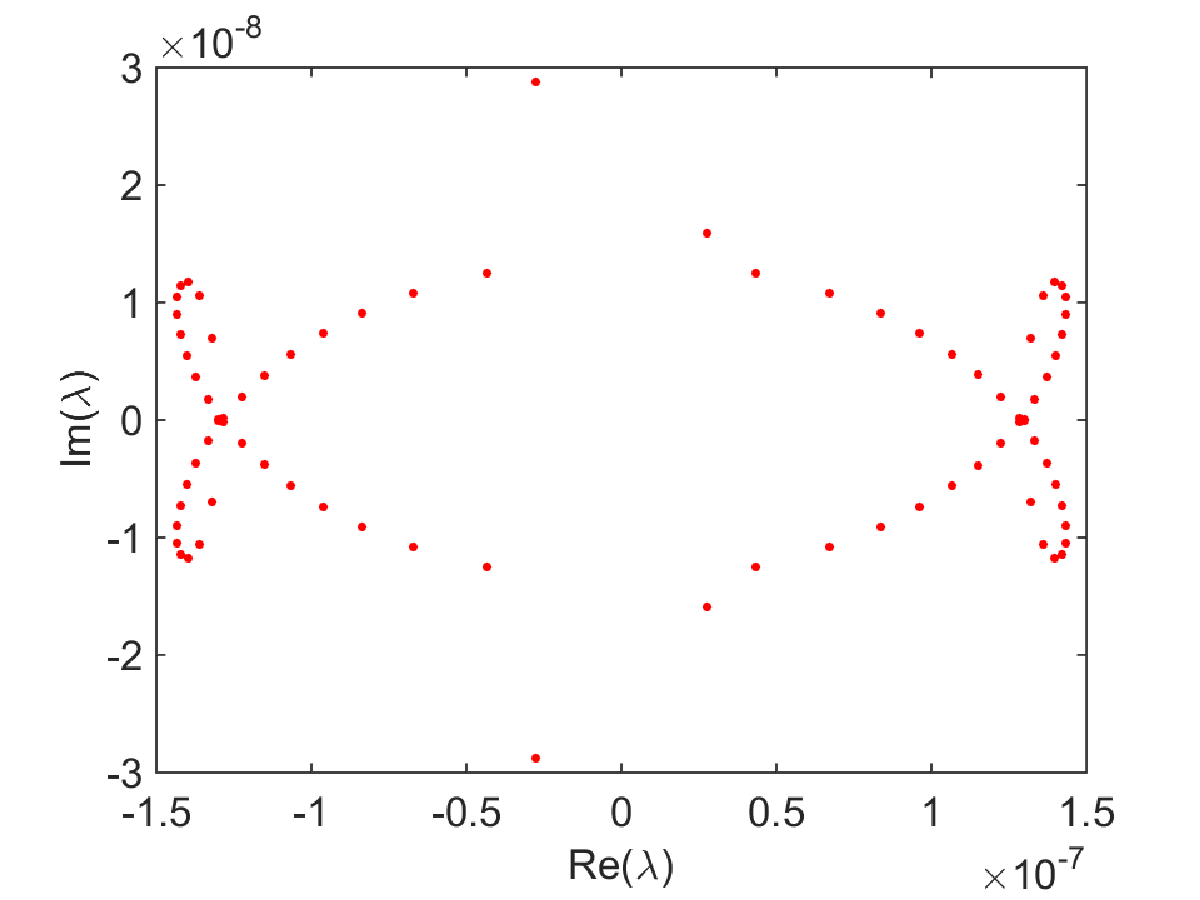}
\includegraphics[width=0.32\textwidth]{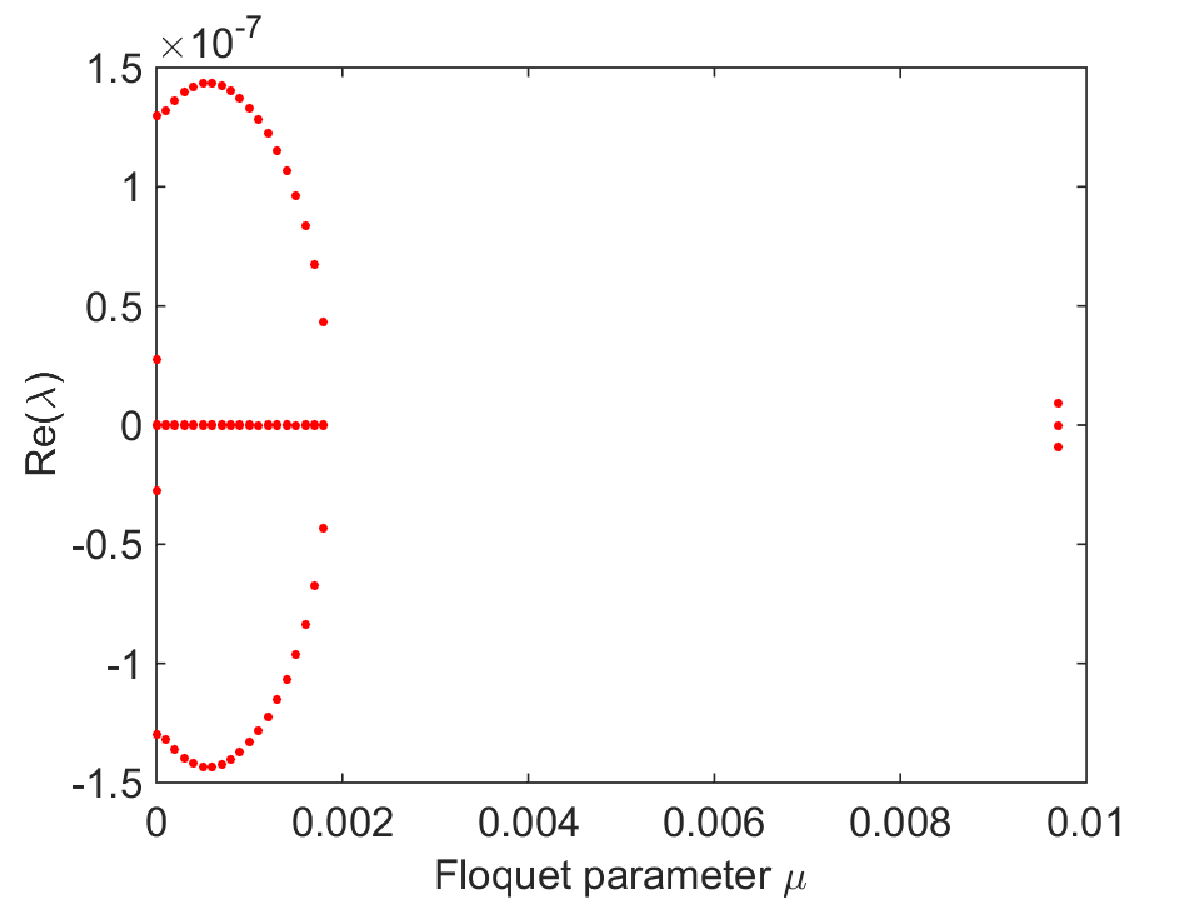}\\
\includegraphics[width=0.32\textwidth]{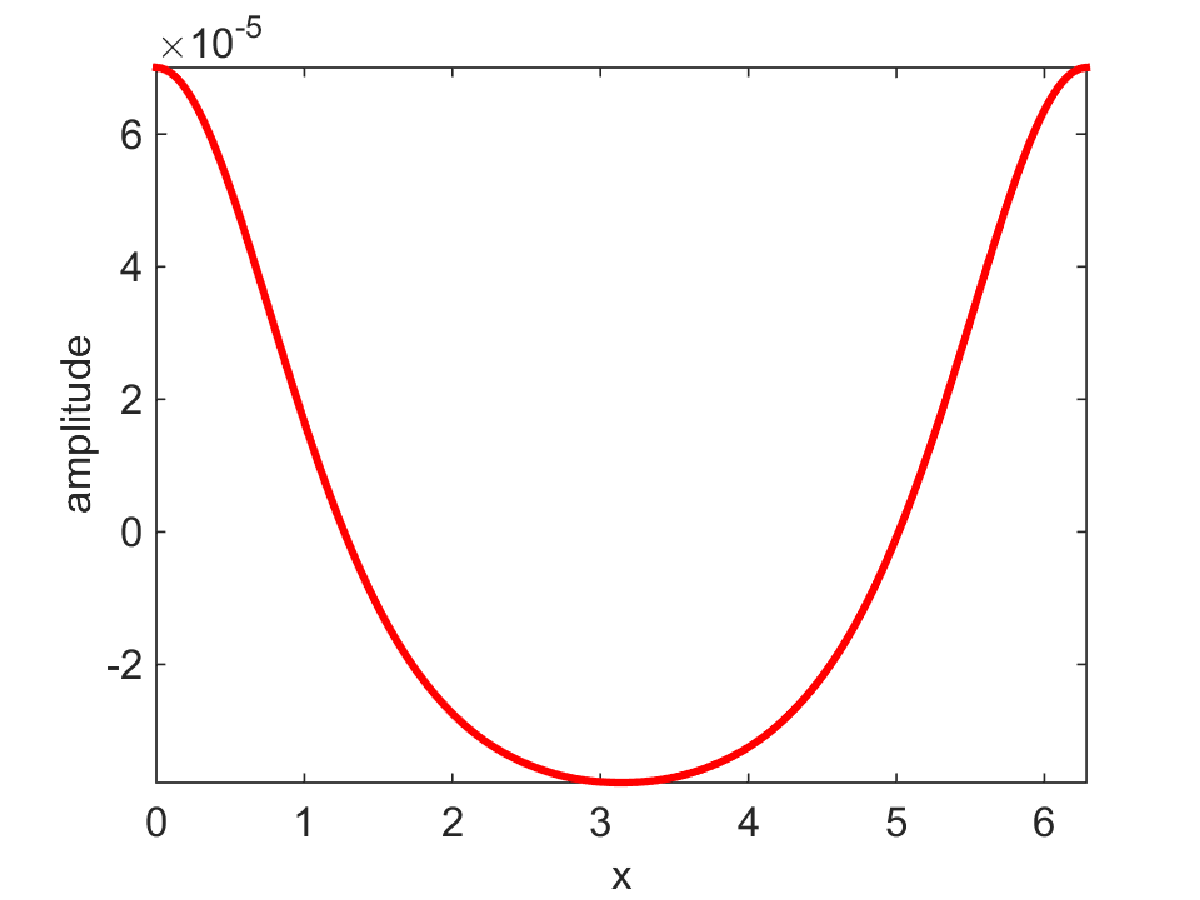}
\includegraphics[width=0.32\textwidth]{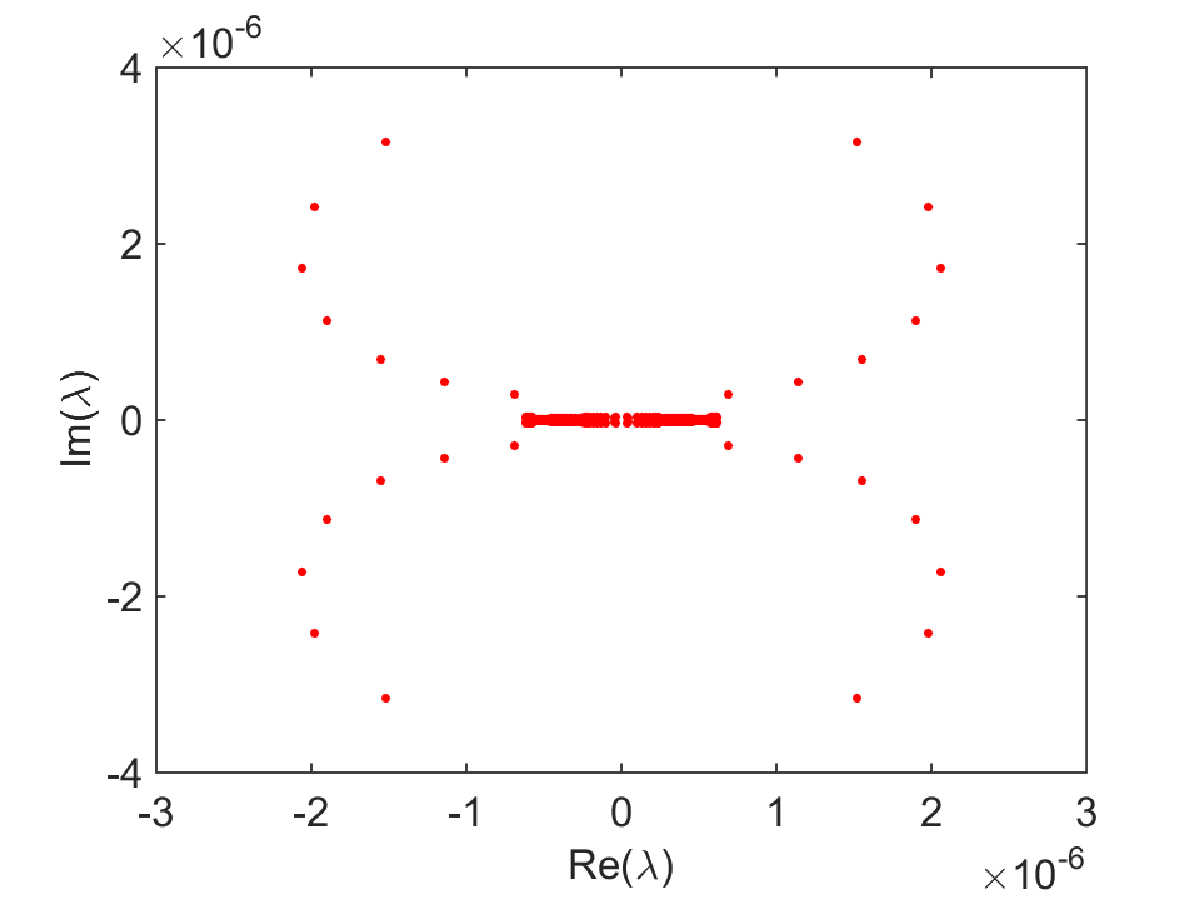}
\includegraphics[width=0.32\textwidth]{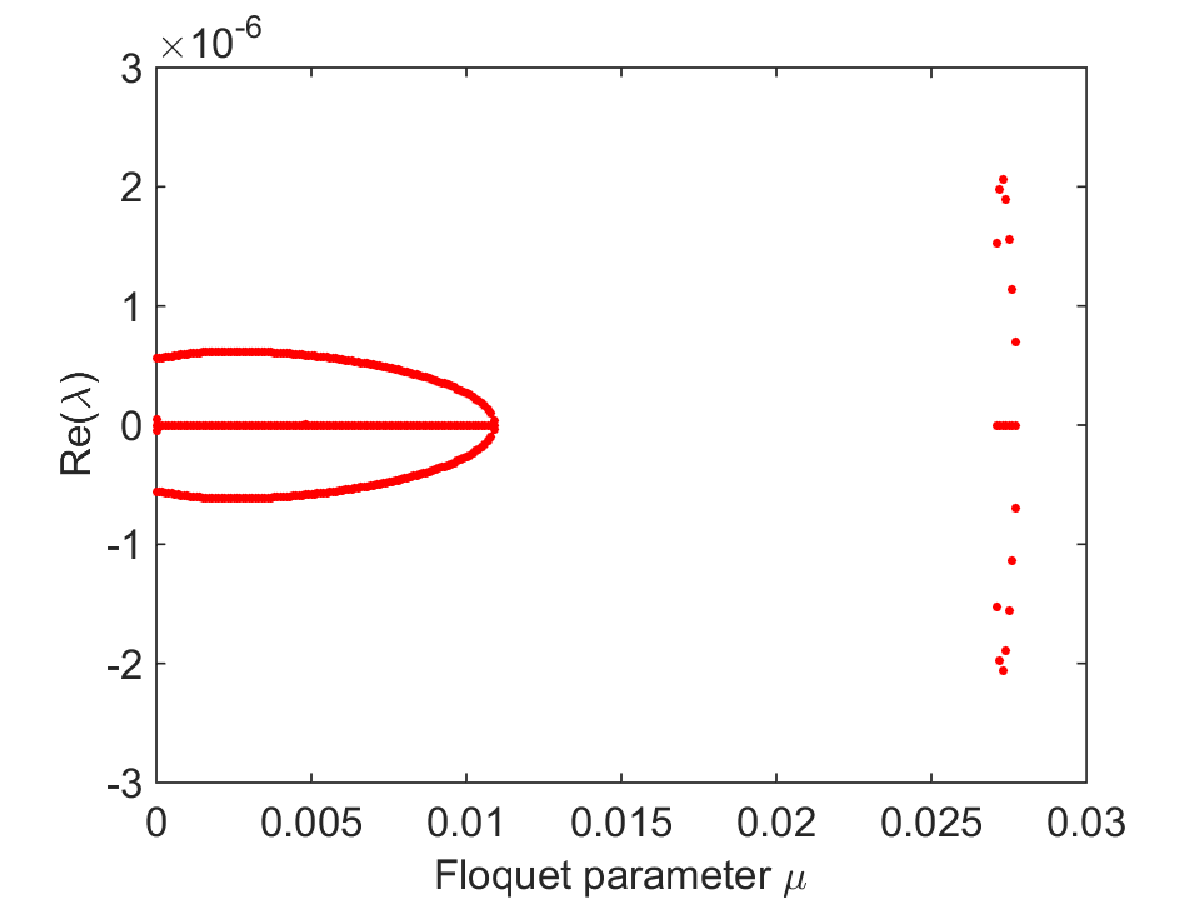}\\
\includegraphics[width=0.32\textwidth]{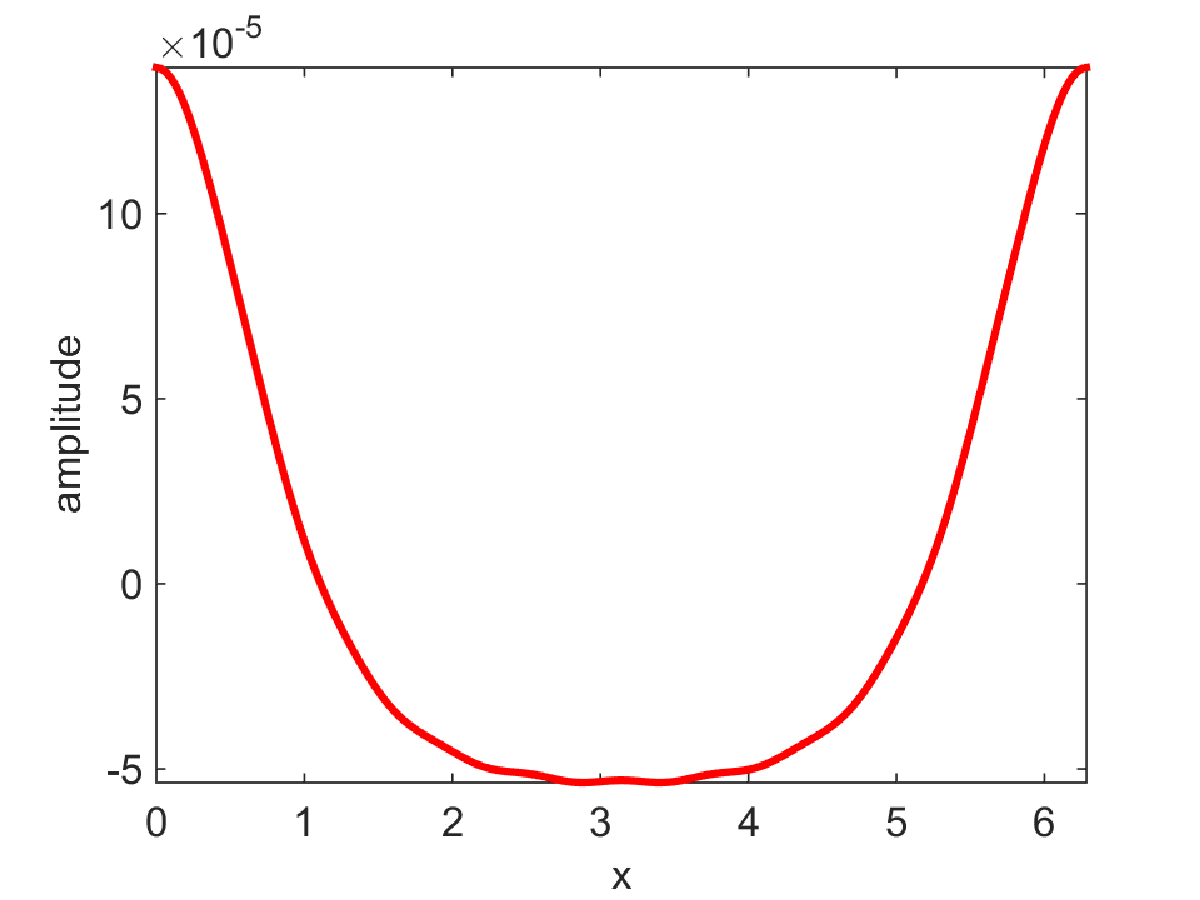}
\includegraphics[width=0.32\textwidth]{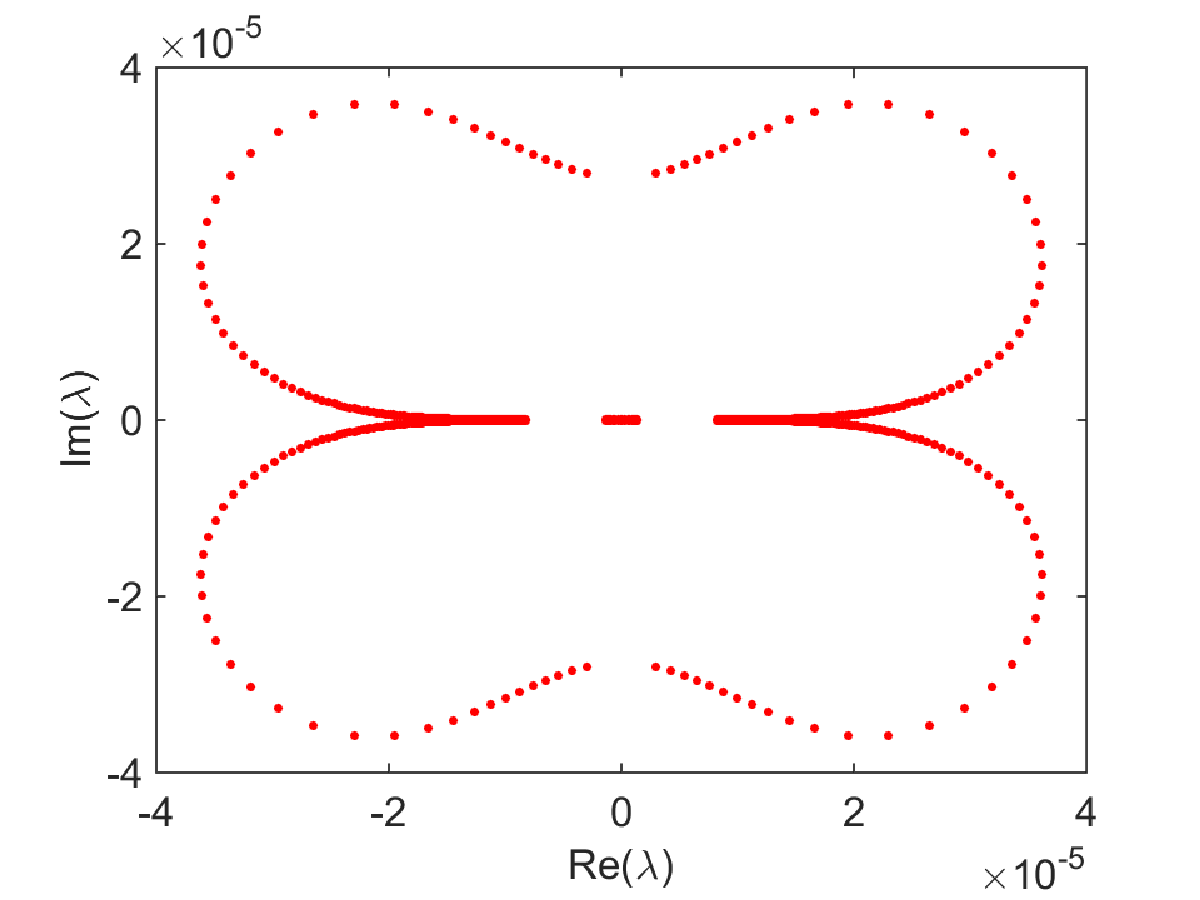}
\includegraphics[width=0.32\textwidth]{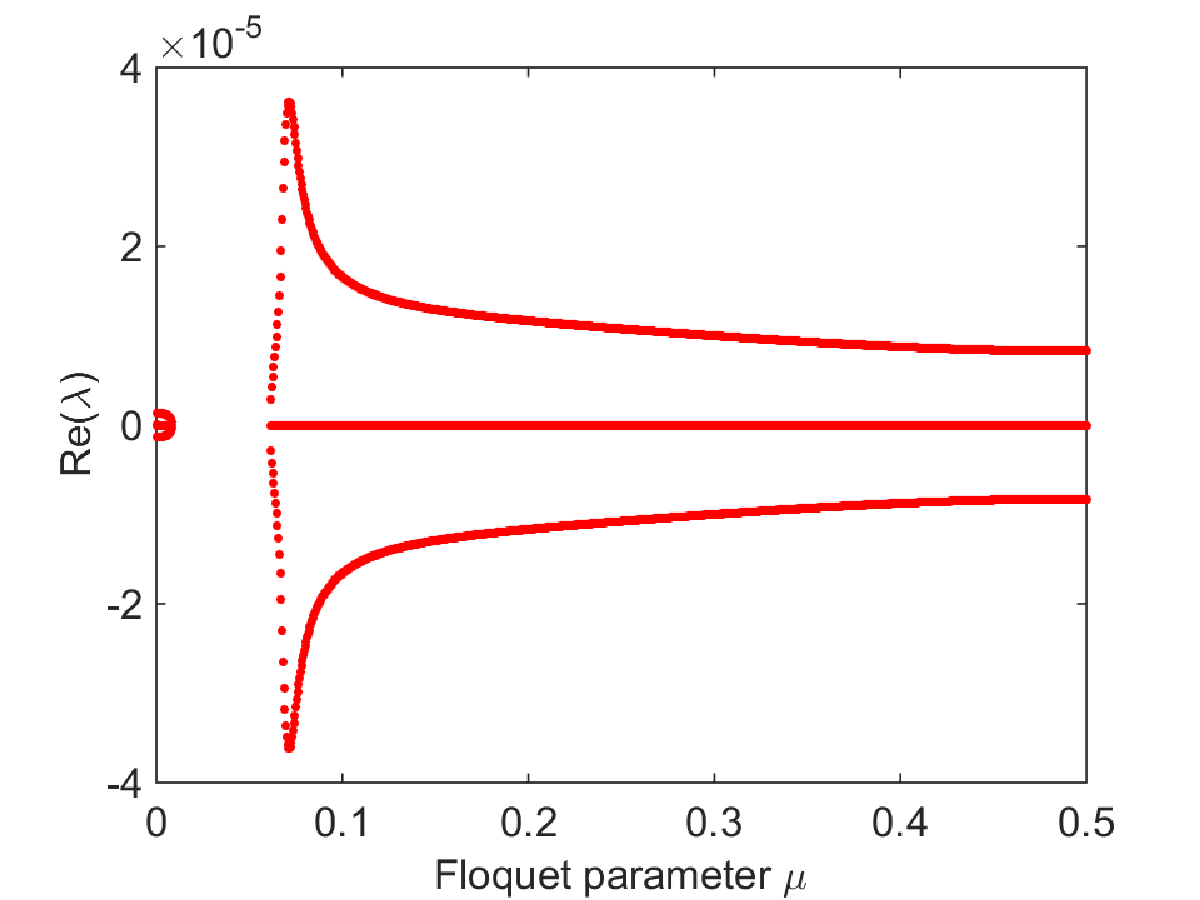}\\
\includegraphics[width=0.32\textwidth]{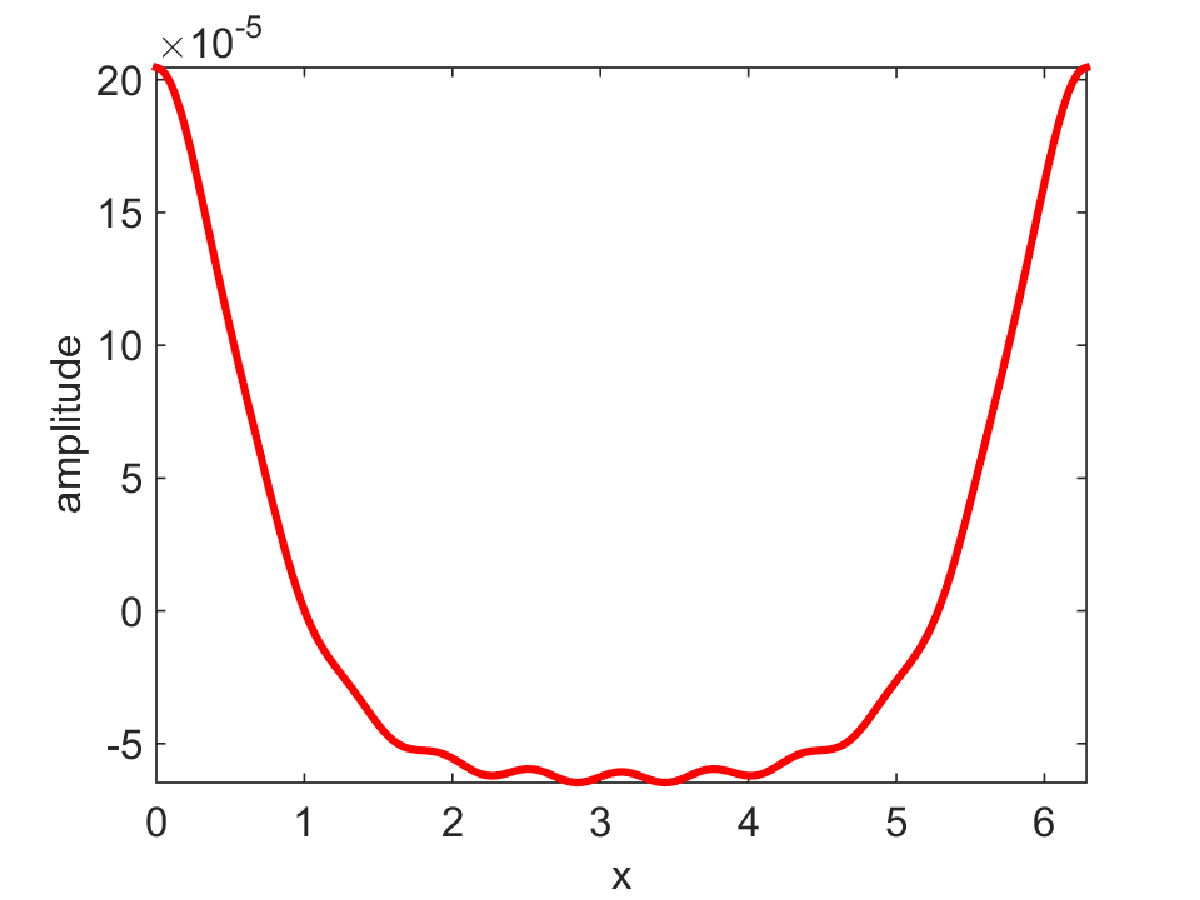}
\includegraphics[width=0.32\textwidth]{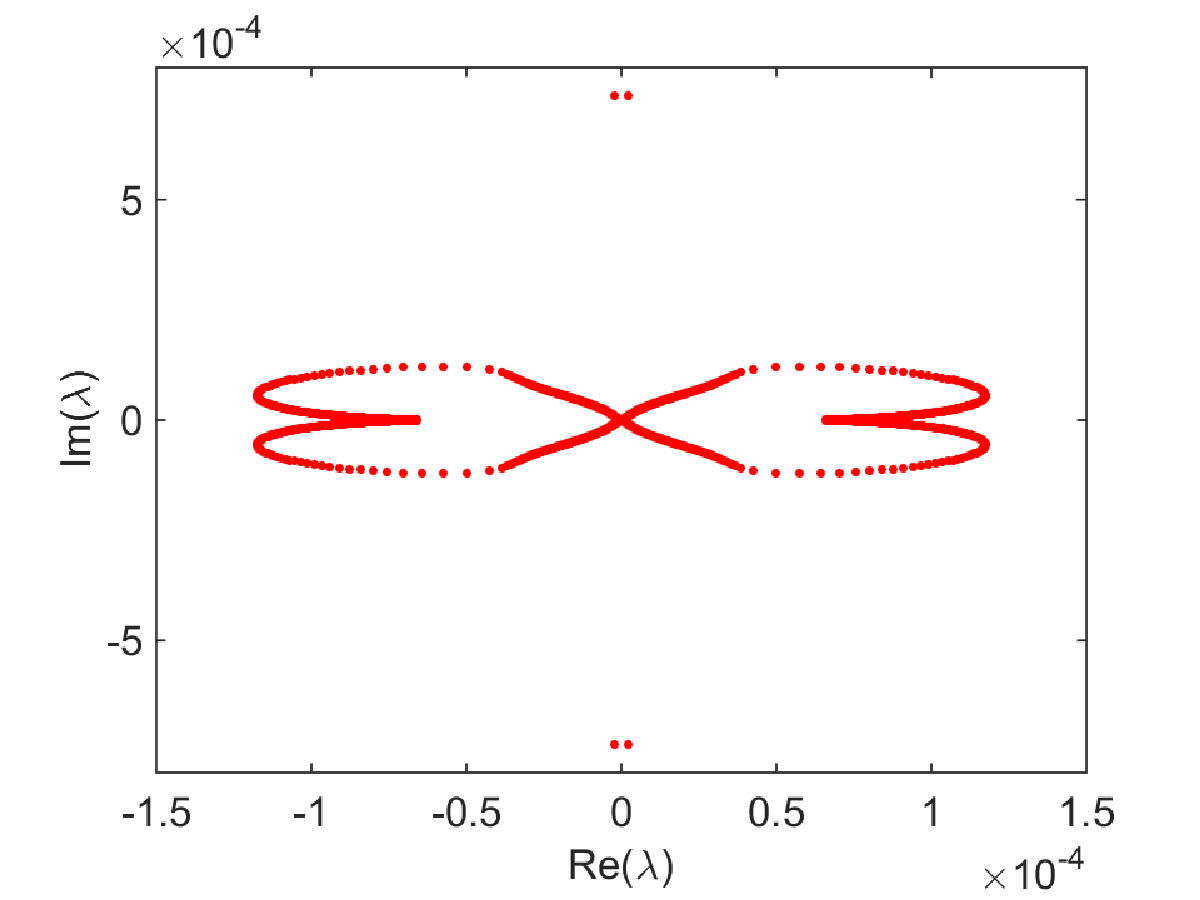}
\includegraphics[width=0.32\textwidth]{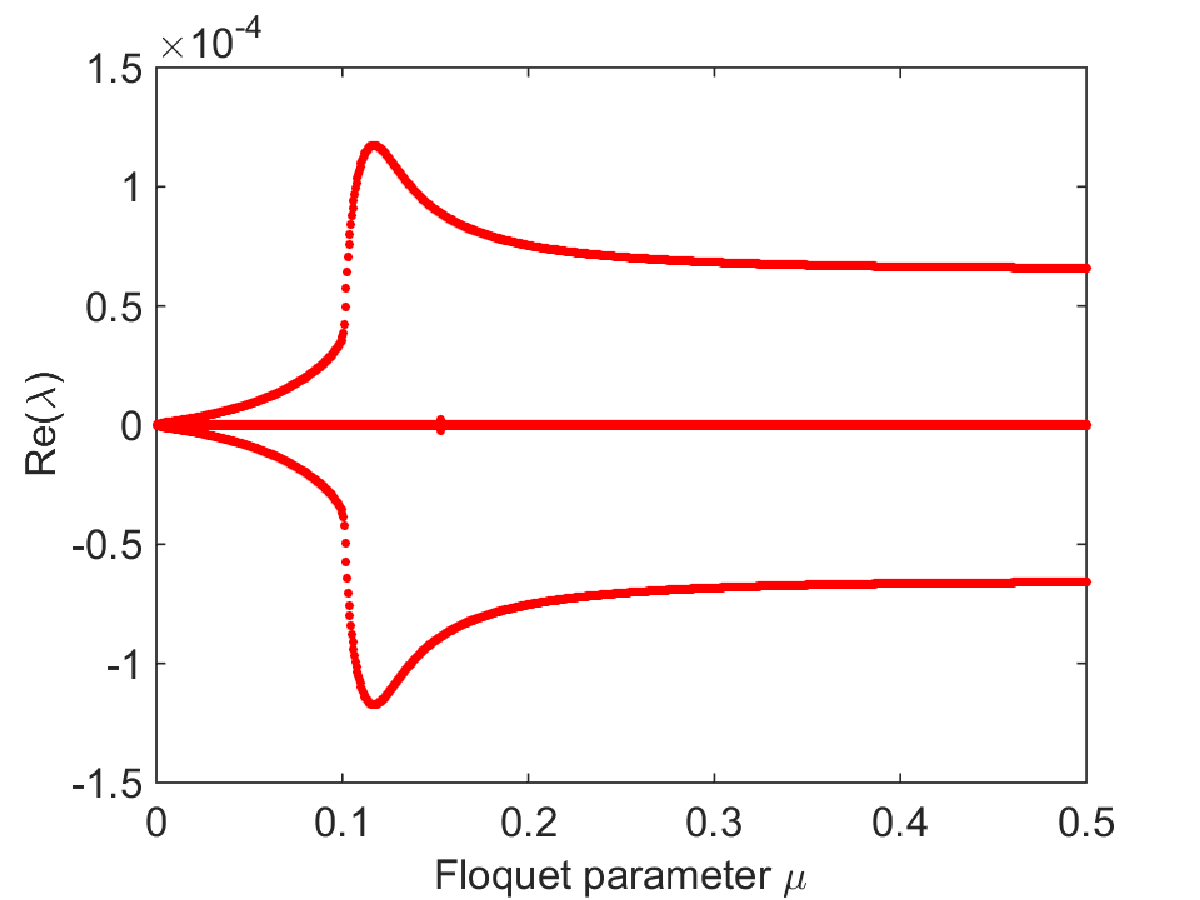}\\
\includegraphics[width=0.32\textwidth]{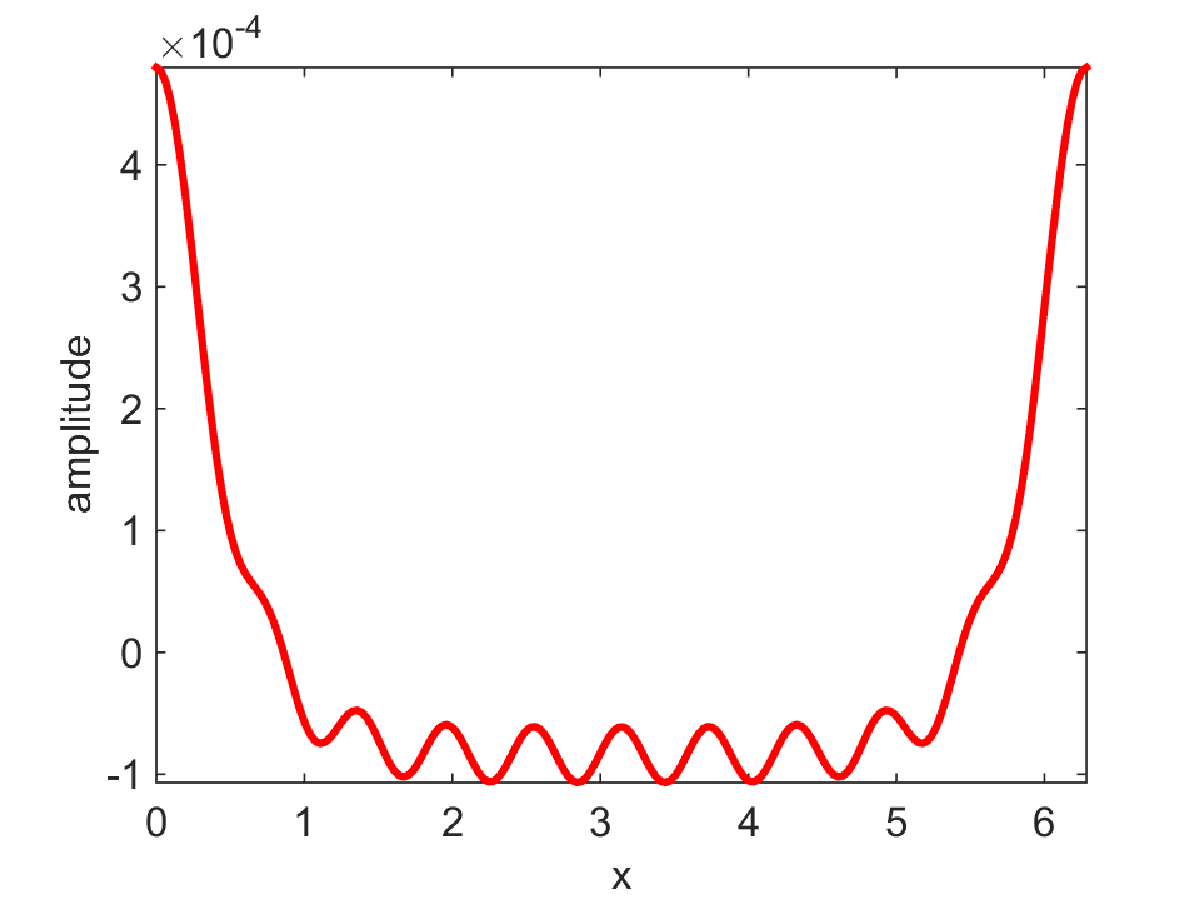}
\includegraphics[width=0.32\textwidth]{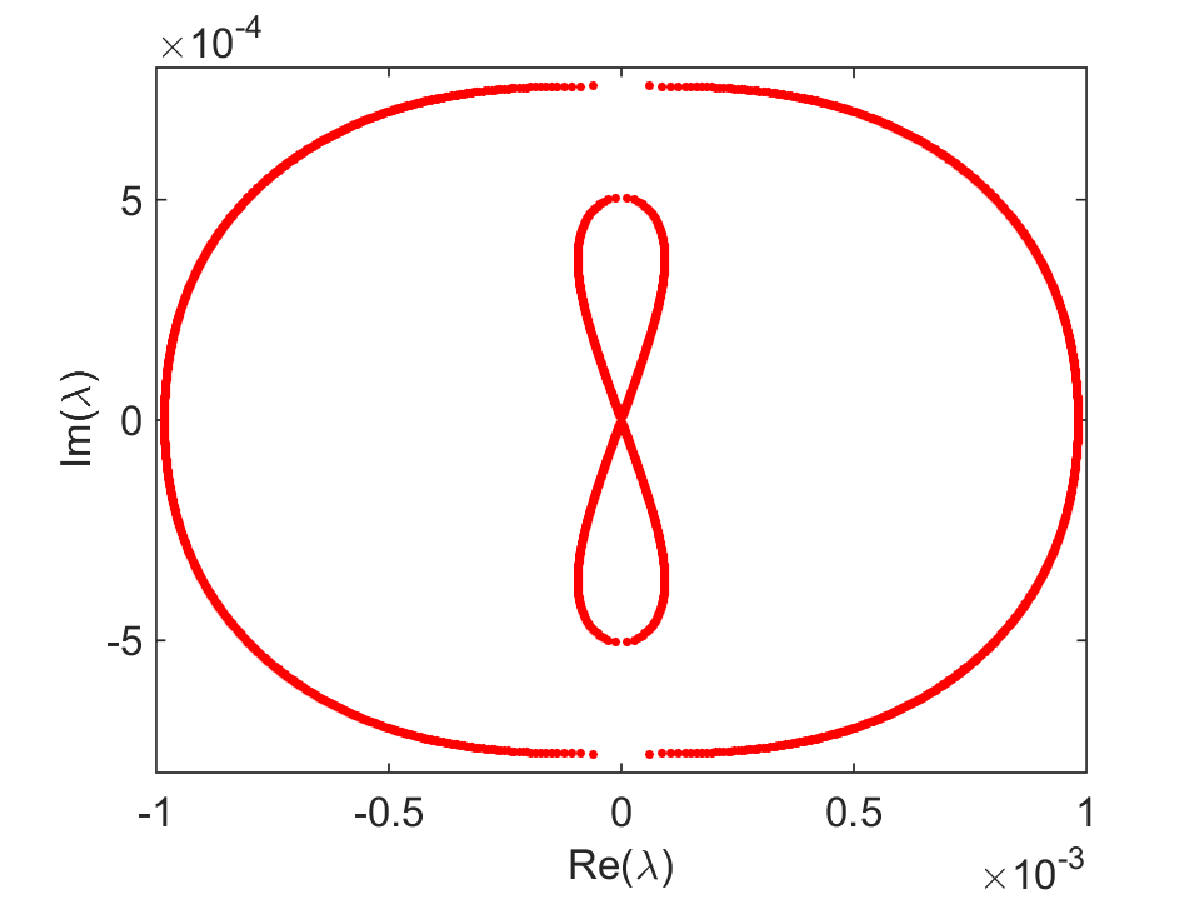}
\includegraphics[width=0.32\textwidth]{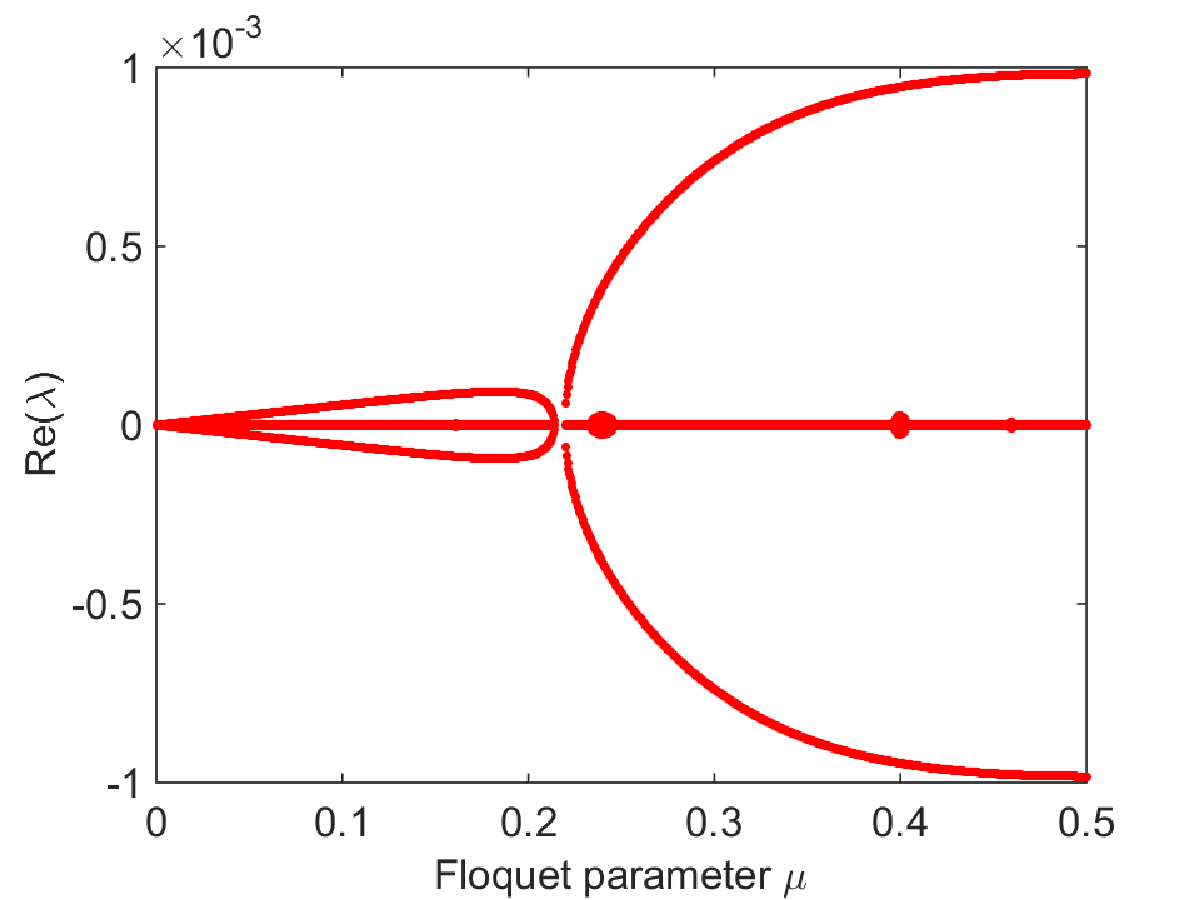}
\caption{Sequence of panels showing waves that increase in amplitude towards the bottom row showing how the modulational intsability arises outside the regime described by asymptotics. Solutions (left) with $K=10$ and $h=0.05$, complex eigenvalue plane (middle) and growth rate versus FLoquet parameter (right). \label{fig:stabK10}}
\end{figure}

\begin{figure}
\includegraphics[width=0.32\textwidth]{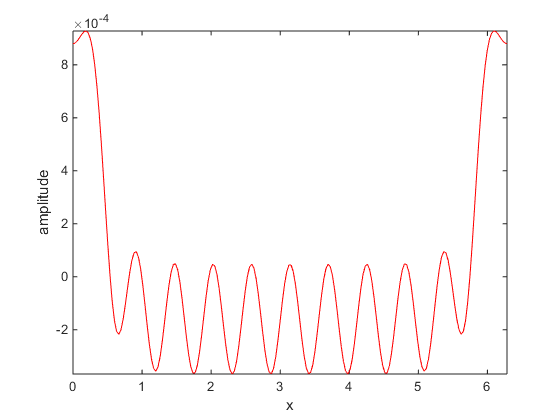}
\includegraphics[width=0.32\textwidth]{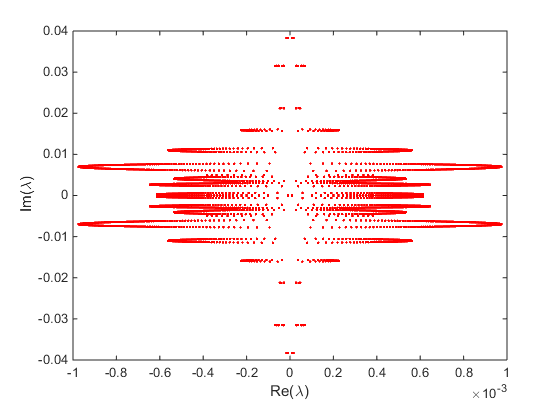}
\includegraphics[width=0.32\textwidth]{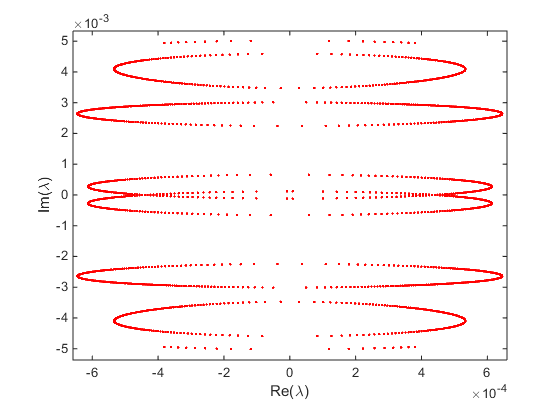}
\caption{Large amplitude resonant solution (left panel) with $K=10$ and $h=0.05$ in the regime where the complex eigenvalue plane (middle panel) shows only high frequency bubble instabilities even near the origin (right panel). \label{fig:K10HighFreqInstab}}
\end{figure}

The complete stability results using Hill's method are shown in Figure \ref{fig:stabK10}. We see that for a small amplitude solutions, there are instabilities near the origin as shown in the top row of the figure. As we increase the amplitude of the solution, a modulational instability arises. The very bottom row shows that high frequency instabilities coexist with a modulational instability for a resonant solution. If we increase the amplitude of the solution even further, we obtain only high frequency instabilities as shown in Figure \ref{fig:K10HighFreqInstab}. 

\begin{figure}
\begin{center}
\begin{tikzpicture}
\node[inner sep=0] at (0,0) {\includegraphics[width=0.32\textwidth]{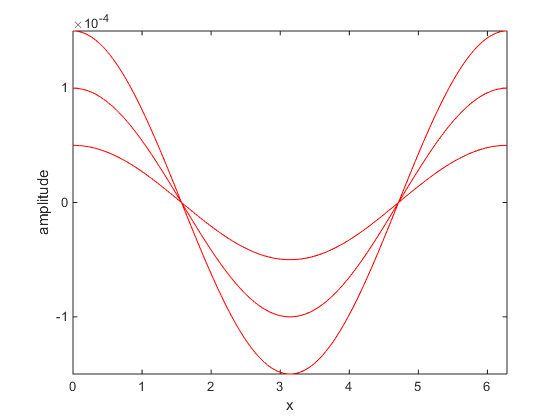}};
\node[inner sep=0] at (5,0) 
{\includegraphics[width=0.32\textwidth]{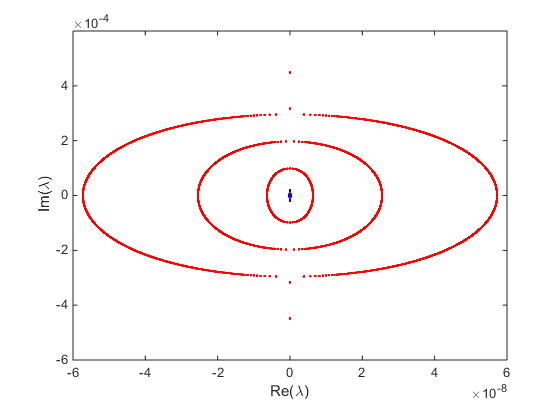}};
\node[inner sep=0] at (10,0) 
{\includegraphics[width=0.32\textwidth]{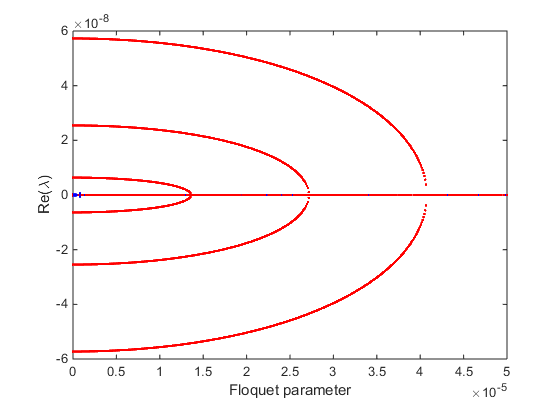}};
    \draw (2.35,0.6) circle (0.165cm);
    \node at (2.35,0.6) {1};
    \draw (2.35,1.1) circle (0.165cm);
    \node at (2.35,1.1) {2};
    \draw (2.35,1.6) circle (0.165cm);
    \node at (2.35,1.6) {3};
    \draw (5.3,0.1) circle (0.165cm);
    \node at (5.3,0.1) {1};
    \draw (5.95,0.1) circle (0.165cm);
    \node at (5.95,0.1) {2};
    \draw (7.0,0.1) circle (0.165cm);
    \node at (7.0,0.1) {3};
    \draw (9.2,0.1) circle (0.165cm);
    \node at (9.2,0.1) {1};
    \draw (10.25,0.1) circle (0.165cm);
    \node at (10.25,0.1) {2};
    \draw (11.3,0.1) circle (0.165cm);
    \node at (11.3,0.1) {3};
\end{tikzpicture}
\end{center}
\caption{Solutions for $D=25$ and $h=\infty$ (left panel) that exhibit bubble instabilities (middle panel) for small Floquet parameters (right panel). \label{fig:InstabD25}}
\end{figure}

We conclude this section by remarking on the very rigid flexural rigidity limit (large $D$) in infinitely deep water. We showed that this regime was asymptotically different for different models in the NLS regime. When analysing the stability of such waves numerically, the instability near the origin is manifested as a bubble instead of what is expected for the modulational instability as shown in Figure \ref{fig:InstabD25}.  There are several explanations for this, stemming from the same phenomenon. This large $D$ limit is near the asymptote as shown in Figure \ref{fig:resonance}. Numerically, this will imply that the first Fourier mode will grow much faster than the others. This is seen when contrasting Figure \ref{fig:deepSolnsD0_1} for $D=0.1$ with Figure \ref{fig:deepSolnsD25} for $D=25$, mainly the large separation in the magnitude of the first and second Fourier modes in the bottom right of the figure for $D=25$. In turn, this implies the assumption of the dependence of modes on a small parameter differs from the one use for deriving NLS. Also as mentioned before, increasing $D$ forces eigenvalues to collide closer and closer to the origin as shown in Figure \ref{fig:eCollisions}. This results in what is seen as a bubble instability in Figure \ref{fig:InstabD25}. The numerical computations are for nonlinear solutions, whereas the asymptotics assumed a linearisation. This means the modulational instability was not seen numerically in this regime for the nonlinear (Toland) model for flexural-gravity waves in its usual form.

\section{Conclusion}\label{sec:conclusion}
Using the AFM reformulation but with two different models describing flexural-gravity waves, we were able derive the local and nonlocal equations for travelling waves under a sheet of ice. By focussing on the travelling wave solutions, we narrowed this down to one equation which was then solved numerically in Fourier space. Assuming an infinite depth, we derived the nonlinear Schr\"{o}dinger equation describing the modulational instability asymptotically. The focussing and defocussing regimes derived using this reformulation with correspond to those seen in \cite{MW13}, but with a different non-dimensionalisation. We showed that the two different models for ice exhibit different stability properties for a large parameter of flexural rigidity $D$ within the NLS regime. We have also confirmed this numerically by first computing solutions to  the Euler's equations and then analysing their stability via the Fourier-Floquet-Hill method. In addition, we examined the resonant regime of the solutions obtained by setting the flexural rigidity parameter such that we obtain a different number of ripples in the wave profile. This effect also does not depend on the model for the  ice, but it is rather apparent in the linear dispersion relation.  When considering high frequency instabilities for waves in finite depth, we showed that these occur in a similar way for either model for the  ice.

\section{Acknowledgements}
This work was supported by EP/J019305/1 for E.P., EP/J019321/1 for P.M. EP/J019569/1 for J.-M.V.-B and O.T.. We would like to thank John Carter for very useful discussions.


\newpage

\bibliographystyle{plain}

\bibliography{biblio}

\end{document}